\documentclass[11pt, oneside]{article}
\usepackage[utf8]{inputenc}
\usepackage{multicol,graphicx,color}
\usepackage{pslatex}
\usepackage{authblk}
\usepackage{amsthm}
\usepackage{amsmath}
\usepackage{amssymb}
\usepackage{latexsym}
\usepackage{lscape}
\usepackage{epsfig}
\usepackage{pstricks}
\usepackage{amsfonts}
\usepackage{mathrsfs}
\usepackage{enumitem}
\usepackage[
  hmarginratio={1:1},     
  vmarginratio={1:1},     
  textwidth=15cm,        
  textheight=21cm,
  heightrounded,          
]{geometry}

\usepackage{graphicx,color}
\usepackage[colorlinks]{hyperref}
\hypersetup{linkcolor=blue,citecolor=blue,filecolor=black,urlcolor=blue}

\usepackage{csquotes}
\usepackage[normalem]{ulem}

\theoremstyle{plain}
\newtheorem{theorem}{Theorem}
\newtheorem{corollary}[theorem]{Corollary}
\newtheorem{lemma}[theorem]{Lemma}

\theoremstyle{definition}
\newtheorem{definition}[theorem]{Definition}
\newtheorem{remark}[theorem]{Remark}

\numberwithin{equation}{section}
\numberwithin{theorem}{section}

\newcommand{\R}{\mathbb{R}}
\newcommand{\N}{\mathbb{N}}
\newcommand{\eps}{\varepsilon}
\newcommand{\dist}{{\rm dist}}

\newcommand{\Zd}{\mathbb{Z}_2}

\newcommand{\M}{S_{\mu}}

\newcommand{\va}{\varepsilon}
\DeclareMathOperator{\sinc}{sinc}
\newcommand{\inv}{S} 

\author{Jack Borthwick\footnote{jack.borthwick@math.cnrs.fr}}  \affil{Laboratoire de Math\'ematiques (CNRS UMR 6623), Universit\'e de Bourgogne Franche-Comt\'e,
Besan\c con 25030, France}

\author{Xiaojun Chang \footnote{chang(x)j100@nenu.edu.cn}} \affil{School of Mathematics and Statistics \& Center for Mathematics and Interdisciplinary Sciences,
 Northeast Normal University, Changchun 130024, Jilin,
PR China}

\author{Louis Jeanjean\footnote{louis.jeanjean@univ-fcomte.fr}}  \affil{Laboratoire de Math\'ematiques (CNRS UMR 6623), Universit\'e de Bourgogne Franche-Comt\'e,
Besan\c con 25030, France}

\author{Nicola Soave\footnote{nicola.soave@unito.it}}  \affil{Dipartimento di Matematica ``Giuseppe Peano", Universit\`a degli Studi di Torino, Via Carlo Alberto 10, 10123, Torino, Italy}

\title{Bounded Palais-Smale sequences with Morse type information for some constrained functionals}
\date{}

\begin{document}

\maketitle

\begin{abstract}
\noindent In this paper, we study, for functionals having a minimax geometry on a constraint, the existence of bounded Palais-Smale sequences carrying  Morse index type information.
\end{abstract}

\medskip

{\small \noindent \text{Keywords:} Critical point theory; Palais-Smale sequences; Morse type information; Variational methods.}

\medskip

{\small \noindent \text{Acknowledgements:}
J. Borthwick gratefully acknowledges that part of this work was supported by the French “Investissements d’Avenir” program, project ISITE-BFC (contract ANR- 15-IDEX-0003).\\
X. J. Chang is partially supported by NSFC (11971095).\\
Part of this work has been carried out in the framework of the Project NQR (ANR-23-CE40-0005-01), funded by the French National Research Agency (ANR). L. Jeanjean thank the ANR for its support.\\
N. Soave is partially supported by the PRIN 2022 project 2022R537CS \emph{$NO^3$
- Nodal Optimization, NOnlinear
elliptic equations, NOnlocal geometric problems, with a focus on regularity} (European Union - Next Generation EU), and by the INdAM - GNAMPA
group.}

\section{Introduction and main results}\label{intro}
In recent years, the search of critical points for constrained functionals having a minimax geometry restricted to a constraint, typically a mountain pass geometry, has been an active direction of research. A core motivation for this is the search for solutions having a prescribed $L^2$-norm to some stationary nonlinear Schr\"odinger equations whose nonlinearity is so-called {\it mass supercritical}.  As a simple model one may consider
\begin{equation}\label{eq:model}
- \Delta u + \lambda u + V(x) u = f(u), \quad u \in H^1(\mathbb{R}^N), 
\end{equation}
where $V(x)$ is a given potential and $f(u)$ is modeled on $f(u) = |u|^{p-2}u$ where  $ 2 + \frac{N}{4} < p < 2^*$ with $2^* = \frac{2N}{N-2}$ and $2^* = + \infty$ if $N=1,2$. A critical point of the associated Energy functional 
$$E(u) = \frac{1}{2} \int_{\R^N}\Big(|\nabla u|^2 + V(x) |u|^2 \Big)\, dx - \frac{1}{p} \int_{\R^N} |u|^{p} \, dx$$
 restricted to the constraint
$$\M = \{u \in H^1(\R^N) \,|\, \int_{\R^N}|u|^2 \, dx = \mu\}$$
for some $\mu>0$ 
is a solution to Equation \eqref{eq:model} with norm $\sqrt{\mu}$ ; the $\lambda \in \R$ arising then as a Lagrange parameter.
Assuming, for simplicity, that  $f(u) = |u|^{p-2}u$ and since $p > 2 + \frac{4}{N}$, the functional $E$  is unbounded from below on $\M$, for any $\mu>0$ but it enjoys a mountain-pass geometry. Namely, there 
exist $w_1, w_2 \in \M$  such that, setting
\[
\Gamma= \left\{ \gamma \in C([0,1],\M) \,|\, \ \gamma(0) = w_1, \quad \gamma(1) = w_2\right\},
\]
we have
\begin{equation}\label{mp-introduction}
c:= \inf_{\gamma \in \Gamma} \ \max_{ t \in [0,1]} E(\gamma(t)) > \max\{E(w_1), E(w_2)\}. 
\end{equation} 
It is a direct consequence of Ekeland's variational principle that a functional with a mountain pass geometry admits Palais-Smale sequences at the mountain pass level $c$. If one can prove that one of these sequences converges then one obtains a critical point at the mountain pass level. However, for constrained functionals, proving the boundedness of such sequences is already a major difficulty. Furthermore, boundedness is of course not sufficient to obtain convergence, and one must also overcome other delicate steps. In particular, using the terminology introduced by P.L. Lions~\cite{Lions84-1,Lions84-2}, one must show that such sequences are not vanishing.
 \medskip

Starting with \cite{Je97}, several approaches have been developed to overcome these issues but ultimately they all seem to rely on the fact that critical points  satisfy a {\it natural} constraint on which the functional can be shown to be coercive and for which, additionally,  sequences {\it close} to the constraint do not vanish, see for example \cite{BaSo17,BeJeLu13}.  In practice, this constraint is often provided by a Pohozaev type identity, or some connected identity, see for example \cite{BaSo17,BeJeLu13, BiMe21, MeSc22, So1} in that direction. Nevertheless, to prove useful, this Pohozaev type constraint must have a simple expression which, in most cases, requires that the underlying functional enjoys some scaling properties. Typically when dealing with problems of the type given by Equation \eqref{eq:model}, these equations must be autonomous or some smallness assumptions on the potential $V(x)$ or its derivatives, which somehow guarantee that the problem keeps some key features of the autonomous case, must be assumed. As a consequence, general non-autonomous equations and, more generally, problems where scaling is not possible -- such as $L^2$-prescribed solutions for mass supercritical equations of type \eqref{eq:model} set on bounded domains or on graphs -- remain essentially unexplored. See nevertheless 
\cite{AlSh22, BaMoRiVe21, MoRiVe, NoTaVe14} for interesting contributions in that direction. \medskip

A central aim of our paper is to present a general abstract tool, Theorem \ref{thm: 0-application}, to deal with cases where a Pohozaev type identity either does not exist or does not provide useful information.  Roughly speaking, we shall show that one can expect generically,  for a functional with a mountain-pass geometry on a constraint, to find a bounded Palais-Smale sequence which carries Morse type information. Along with that of boundedness, this second-order information will be key in the proof that the sequence converges. In fact, our Theorem \ref{thm: 0-application} has already proved decisive in the papers \cite{BoChJeSo:2022,ChJeSo:2022} in which we study mass prescribed solutions for a mass-superlinear problem set on a graph. The fact that the underlying equation is set on a graph prevents the existence of a useful Pohozaev identity. Actually, Theorem \ref{thm: 0-application} will be derived as a corollary of a more general result, Theorem \ref{thm: application}, in which more general minimax geometries than the mountain pass one can be handled.    \medskip

\subsection*{Setting and statement of Theorem~\ref{thm: 0-application} and Theorem~\ref{thm: application}}
We now introduce the setting in which our main results will be stated. This setting is borrowed from \cite[Section 8]{BerestyckiLions1983-2}. \medskip\\
Let $(E,\langle \cdot, \cdot \rangle)$ and $(H,(\cdot,\cdot))$ be two \emph{infinite-dimensional} Hilbert spaces and assume that:
\[ E\hookrightarrow H \hookrightarrow E',\]
with continuous\footnote{In the applications, these injections will also be dense but we have not used this property.} injections.
 For simplicity, we assume that the continuous injection $E\hookrightarrow H$ has norm at most $1$ and identify $E$ with its image in $H$. We also introduce:  \[ \begin{cases} \|u\|^2=\langle u,u \rangle, \\ |u|^2=(u,u),\end{cases}\quad u\in E.\]
Define for $\mu \in (0,+\infty)$: \begin{equation}\label{def-M} \M= \{ u \in E \,|\, |u|^2=\mu \}. \end{equation}
Note that $\M$ is a submanifold of $E$ of codimension $1$ and that its tangent space at a given point $u \in \M$ can be considered as the closed codimension $1$ subspace of $E$ given by:
\[  T_u \M = \{v \in E \,|\, (u,v) =0 \}.\]

In order to state Theorem \ref{thm: 0-application} we need some definitions. We denote by $\|\cdot\|_*$ and $\|\cdot\|_{**}$, respectively, the operator norm of $\mathcal{L}(E,R)$ and of $\mathcal{L}(E,\mathcal{L}(E,R))$.

\begin{definition}\label{Holder continuous}
Let $\phi : E \rightarrow \mathbb{R}$ be a $C^2$-functional on $E$ and $\alpha \in (0,1]$; we shall say that $\phi'$ and $\phi''$ are:
\begin{enumerate}[label=\textbf{HOL \arabic*}] \item:\label{hol1} \textbf{(globally) $\alpha$-Hölder continuous} if there is $M>0$ such that for any $u_1,u_2 \in E$,
\[||\phi'(u_1)-\phi'(u_2)||_*\leq M ||u_2-u_1||^{\alpha}, \quad ||\phi''(u_1)-\phi''(u_2)||_{**}\leq M ||u_1-u_2||^\alpha, \]
\item:\label{hol2} \textbf{$\alpha$-Hölder continuous on bounded sets} if for any $R>0$ one can find $M=M(R)>0$ such that for any $u_1,u_2\in B(0,R)$:
\begin{equation}\label{Holder}
||\phi'(u_1)-\phi'(u_2)||_*\leq M ||u_2-u_1||^{\alpha}, \quad ||\phi''(u_1)-\phi''(u_2)||_{**}\leq M||u_1-u_2||^\alpha,
\end{equation}
\item:\label{hol3} \textbf{locally $\alpha$-Hölder continuous} if for any $u\in E$ one can find an open neighbourhood $U$ of $u$ such that the restriction of $\phi$ to $U$, $\phi|_U$, satisfies $\ref{hol1}$.
\end{enumerate}
\end{definition}

\begin{remark}
Clearly:  \ref{hol1} $\Rightarrow$ \ref{hol2} $\Rightarrow$ \ref{hol3}. 
 \end{remark}
\begin{definition}\label{def D}
Let $\phi$ be a $C^2$-functional on $E$, for any $u\in E\setminus\{0\}$ define the continuous bilinear map:
\[ D^2\phi(u)=\phi''(u) -\frac{\phi'(u)\cdot u}{|u|^2}(\cdot,\cdot).  \]
\end{definition}


\noindent
The geometric relevance of Definition~\ref{def D} will be explained in Section~\ref{relevanceD2}.

\begin{definition}\label{def: app morse}
Let $\phi$ be a $C^2$-functional on $E$, for any $u\in \M$ and $\theta >0$, we define
 the \emph{approximate Morse index} by
\[
\tilde m_\theta(u) = \sup \left\{\dim\,L\left|~L \text{ is a subspace of $T_u \M$ such that:~$\forall \varphi \in L \backslash \{0 \}, \,$ } 
D^2\phi(u) [\varphi, \varphi] <-\theta \|\varphi\|^2\right.\right\}.
\]
\end{definition}

Our first main result is the following.

\begin{theorem}\label{thm: 0-application}
Let $I \subset (0,+\infty)$ be an interval and consider a family of $C^2$ functionals $\Phi_\rho: E \to \mathbb{R}$ of the form
\[
\Phi_\rho(u) = A(u) -\rho B(u), \qquad \rho \in I,
\]
where $B(u) \ge 0$ for every $u \in E$, and
\begin{equation}\label{0-hp coer}
\text{either $A(u) \to +\infty$~ or $B(u) \to +\infty$ ~ as $u \in E$ and $\|u\| \to +\infty$.}
\end{equation}
Suppose moreover that $\Phi_\rho'$ and $\Phi_\rho''$ are $\alpha$-Hölder continuous on bounded sets in the sense of Definition~\ref{Holder continuous}~(\ref{hol2}) for some $\alpha \in (0,1]$.
Finally, suppose that there exist $w_1, w_2 \in \M$ (independent of $\rho$) such that, setting
\[
\Gamma= \left\{ \gamma \in C([0,1],\M) \,|\,  \gamma(0) = w_1, \quad \gamma(1) = w_2\right\},
\]
we have
\begin{equation}\label{0-mp geom}
c_\rho:= \inf_{\gamma \in \Gamma} \ \max_{ t \in [0,1]} \Phi_\rho(\gamma(t)) > \max\{\Phi_\rho(w_1), \Phi_\rho(w_2)\}, \quad \rho \in I.
\end{equation}
Then, for almost every $\rho \in I$, there exist sequences $\{u_n\} \subset \M$ and $\zeta_n \to 0^+$ such that, as $n \to + \infty$,
\begin{itemize}
\item[(i)] $\Phi_\rho(u_n) \to c_\rho$;
\item[(ii)] $||\Phi'_\rho|_{\M}(u_n)||_* \to 0$; 
\item[(iii)] $\{u_n\}$ is bounded in $E$;
\item[(iv)] $\tilde m_{\zeta_n}(u_n) \le 1$.
\end{itemize}
\end{theorem}

\begin{remark}\label{Gradient}
For any $u \in \M$ we denote by: $\Phi'_\rho|_{\M}(u)$ the differential at $u$ of the restriction of $\Phi_\rho$ to $\M$, it is therefore a linear map: $T_uS_\mu \rightarrow T_{\Phi_\rho(u)}\R \cong \R$. The notation $||\Phi'_\rho|_{\M}(u)||_*$, introduced before Definition \ref{Holder continuous}, and used in Theorem \ref{thm: 0-application} (ii), is understood to refer to the dual norm induced by the norm of $T_u\M$ which is inherited from $E$. It is known from \cite[Lemma 3]{BerestyckiLions1983-2} that, if a sequence  $\{v_n\} \subset \M$ is bounded, then the following are equivalent:
\begin{itemize}
\item $\|\Phi'_\rho|_{\M}(v_n)\|_* \to 0,\,$ as $n \to + \infty$.
\item $\Phi'_{\rho}(v_n) - \displaystyle \frac{1}{\mu}( \Phi’_\rho(v_n)\cdot v_n)(v_n, \cdot) \to 0 \,$ in  $E'\,$, as $n \to + \infty$.
\end{itemize}
Thus, in particular, Theorem \ref{thm: 0-application}  (ii) implies that
\begin{equation}\label{free-gradient}
 \Phi'_{\rho}(u_n) - \displaystyle \frac{1}{\mu}( \Phi’_\rho(u_n)\cdot u_n)(u_n, \cdot) \to 0 \, \mbox{ in }  E' \mbox{ as } n \to + \infty.
\end{equation}
\end{remark}

\begin{remark}\label{Morse-explicit}
It follows immediately from Theorem \ref{thm: 0-application}  (iv) that if there exists a subspace $W_n \subset T_{u_n}\M$ such that
\[
\Phi''_{\rho}(u_n)[w,w] - \frac{1}{\mu}( \Phi’_\rho(u_n)\cdot u_n)(w,w) < - \zeta_n ||w||^2, \quad \mbox{for all } w \in W_n \backslash \{0 \},
\]
then necessarily $\dim W_n \leq 1$.
\end{remark}

\begin{remark}\label{space-of-functions}
Typically, Theorem \ref{thm: 0-application} will be used with a space $E$ which consists of functions defined on $\R^N$, some subset of $\R^N$, or graphs as in
\cite{BoChJeSo:2022,ChJeSo:2022}.  For instance, to deal with Equation \eqref{eq:model} we shall set $E= H^1(\R^N)$ and $H= L^2(\R^N)$.\\
 In such situations, denoting by $|u|_*$ the modulus of a function $u \in E$, we shall see that it is possible to choose  $\{u_n\} \subset \M$ with the property that $u_n \ge 0$ on 
$E$ if all of the following are satisfied: $u \in E\ \mapsto \ |u|_* \in E$, $w_1, w_2 \ge 0$, the map $u \mapsto |u|_*$ is continuous, and $\Phi_\rho(u) = \Phi_\rho(|u|_*)$. 
\end{remark}

Actually, Theorem \ref{thm: 0-application} will be a consequence of a more general result, Theorem \ref{thm: application} below. To state it we need the following definition.

\begin{definition}\label{C-def: def1.5}
A family $\mathcal{F}$ of subsets of $\M$ is said to be homotopic of dimension at most $d$ with boundary $B$ if there exist a compact subset $D$ of $\mathbb{R}^d$, containing a closed subset $D_0,$ and a continuous function $\eta_0$ from $D_0$ onto $B$ such that
$$
\mathcal{F}= \{ A \subset \M \,|\, A = f(D) \mbox{ for some } f \in C(D;\M) \mbox{ with } f = \eta_0 \mbox{ on } D_0 \}.$$
\end{definition}


\begin{theorem}\label{thm: application}
Let $I \subset (0,+\infty)$ be an interval and consider a family of $C^2$ functionals $\Phi_\rho: E \to \mathbb{R}$ of the form
\[
\Phi_\rho(u) = A(u) -\rho B(u), \qquad \rho \in I,
\]
where $B(u) \ge 0$ for every $u \in E$, and
\begin{equation}\label{hp coer}
\text{either $A(u) \to +\infty$~ or $B(u) \to +\infty$ ~ as $u \in E$ and $\|u\| \to +\infty$.}
\end{equation}
Suppose moreover that $\Phi_\rho'$ and $\Phi_\rho''$ are $\alpha$-Hölder continuous on bounded sets in the sense of Definition~\ref{Holder continuous}~(\ref{hol2}) for some $\alpha \in (0,1]$.

Let $\mathcal{F}$ be a homotopic family of $\M$ of dimension at most $d$ with boundary $B$ (independent of $\rho$) such that 
\begin{equation}\label{mp geom}
c_{\rho}:= \inf_{A \in \mathcal{F}} \max_{u \in A}\Phi_{\rho}(u) > \max_B \Phi_{\rho},\quad \forall \rho \in I.
\end{equation}
Then, for almost every $\rho \in I$, there exist sequences $\{u_n\} \subset \M$ and $\zeta_n \to 0^+$ such that, as $n \to + \infty$,
\begin{itemize}
\item[(i)] $\Phi_\rho(u_n) \to c_\rho$;
\item[(ii)] $||\Phi'_\rho|_{\M}(u_n)||_* \to 0$; 
\item[(iii)] $\{u_n\}$ is bounded in $E$;
\item[(iv)] $\tilde m_{\zeta_n}(u_n) \le d$.
\end{itemize}
\end{theorem}

\begin{remark}\label{linking}
Theorem \ref{thm: application} covers the cases where the family of $C^2$ functionals $\Phi_\rho: E \to \mathbb{R}$ has a uniform mountain-pass, respectively a uniform linking, geometry. Indeed, 
\begin{enumerate}
\item  Suppose that there exist $w_1, w_2 \in \M$ (independent of $\rho$) such that, setting
\[
\Gamma= \left\{ \gamma \in C([0,1];\M) \,|\,  \gamma(0) = w_1, \, \gamma(1) = w_2\right\},
\]
we have
\begin{equation*}
c_\rho:= \inf_{\gamma \in \Gamma} \ \max_{ t \in [0,1]} \Phi_\rho(\gamma(t)) > \max\{\Phi_\rho(w_1), \Phi_\rho(w_2)\}, \quad \rho \in I.
\end{equation*}
Then 
$$
\mathcal{F}= \{A \subset \M \,|\, A= \gamma([0,1]) \mbox{ for some } \gamma \in \Gamma \}
$$
  is a homotopic family of $\M$ of dimension at most $1$ with boundary $\{w_1, w_2\}$, such that (\ref{mp geom}) holds. Thus Theorem \ref{thm: application} indeed extends Theorem \ref{thm: 0-application}.
  \item Assume that 
\begin{itemize}
\item[(i)] $S \subset S_{\mu} $ is a closed subset of $S_{\mu}$.
\item[(ii)] $Q \subset S_\mu$ is such that $Q=\eta(D)$, where $ \eta:D \subset \R^d \to S_\mu$ is continuous and $D \subset \R^d$ is a regular compact set with non-empty interior.
\item[(iii)] denoting by $\partial Q = \eta(\partial D)$ - the \emph{boundary} of $Q$ in $S_\mu$ - we have that $\partial Q$ and $S$ link, namely
$\partial Q \cap S = \emptyset$ and for all $g \in C(\M, \M)$ with $g = Id$ on $\partial Q$, we have $g(Q) \cap S \neq \emptyset$.
\item[(iv)] \begin{equation*}
\sup\limits_{u\in \partial Q}\Phi_\rho(u)< \inf\limits_{u\in S}\Phi_\rho(u), \quad \mbox{for any } \rho\in I.
\end{equation*}
\end{itemize}
Then setting
$$
\mathcal{F}= \{A \subset \M \,|\, A= f(D) \mbox{ for some } f \in C(D;\M) \mbox{ with } f = \eta \mbox{ on } \partial D \}
$$
we have that $\mathcal{F}$ is a homotopic family of $\M$ of dimension at most $d$ with boundary $\partial Q$ such that (\ref{mp geom}) holds.
\end{enumerate}
In fact, the case of a saddle point geometry, see \cite{Rabi1986} for a definition, can also be covered by Theorem \ref{thm: application}.
\end{remark}

%
%

Under the assumptions of Theorem \ref{thm: application}, one can readily observe that the function $\rho \mapsto c_\rho$ is non-increasing. Therefore, its derivative $c_\rho'$ is well defined for almost every $\rho \in I$. Then, using a tool first introduced by Struwe in \cite{Struwe88}, and following the approach of it developed in \cite{J-PRSE1999}, we shall show that the existence of $c_\rho'$ ensures, for such a value of $\rho$, the existence of a sequence of minimising paths in $\Gamma$ whose {\it tops} lie in a given ball. This key observation could lead rather directly to the existence of a bounded Palais-Smale sequence at the level $c_\rho$, but we additionally aim to prove that some of these sequences possess approximate Morse index properties. \medskip

The first author to emphasise the importance of second order conditions, such as ($iv$) in Theorem \ref{thm: application}, in compactness problems was P.-L. Lions \cite{Li87}, see also \cite{BaLi87, BaLi92}.  The main idea in these papers was to use Morse type information on exact critical points. A major contribution in dealing instead with almost critical points, as is the case for a Palais-Smale sequence, is due to Fang and Ghoussoub, see \cite[Theorem 1]{Fang:1992wz}. Therein, they proved (along with other valuable properties) that an unconstrained functional having a mountain-pass geometry possesses, at the mountain-pass level, a Palais-Smale sequence whose elements have {\it approximately} Morse index 1.

Unfortunately, \cite[Theorem 1]{Fang:1992wz} requires that the first and second derivatives of the functional satisfy the Hölder continuity assumption~\ref{hol1} in a neighbourhood of the level set corresponding to the mountain pass value. Actually, inspecting the proof of \cite[Theorem 1]{Fang:1992wz}, one observes that this condition only needs to be required around some {\it dual} set. However, without precise information on the functional, one has to check the full global condition~\ref{hol1} which has limited the possibility to use directly \cite[Theorem 1]{Fang:1992wz}. The ideas of \cite{Fang:1992wz} nevertheless played a crucial role in the works \cite{BeRu} and \cite{LoMaRu} in which, combining the ideas of \cite{Fang:1992wz, J-PRSE1999}, the authors managed to derive unconstrained versions of our Theorem \ref{thm: application}, under the weaker assumption \ref{hol2}.   \medskip

In fact, in \cite{Fang:1994wz}, \cite[Theorem 1]{Fang:1992wz} was extended, directly under the assumption that the H\"older condition~\ref{hol1} holds, to a much larger setting. In particular, in \cite{Fang:1994wz} one can handle general min-max structures, as introduced in Theorem \ref{thm: application} instead of just a mountain pass one. Moreover, it is stated  in \cite{Fang:1994wz} that the result of \cite{Fang:1992wz}  holds for functionals defined on a complete $C^2$- Riemannian manifold modeled on a Hilbert space $E$, see \cite[Theorem 1.7 and Theorem 3.1]{Fang:1994wz} for the precise statements.  However, the proof of these two theorems is only carried out for an (unconstrained) functional defined on a Hilbert space. In~\cite[page 1609]{Fang:1994wz}, the authors indicate that the general case follows by a direct adaptation of the proof given in the Hilbert case. 

Whilst we agree that the arguments used to prove the various lemmata (\cite[Lemma 3.(3-7)]{Fang:1994wz}) that build up to the proof of the main result~\cite[Theorem 3.1]{Fang:1994wz} are local in nature and can be made sense of in a local chart, one cannot substitute the \emph{local} Hölder assumption~\ref{hol3} -- which one could expect to extend in a standard fashion to manifolds -- for~\ref{hol1} directly in the statement of~\cite[Theorem 3.1]{Fang:1994wz}. More specifically, inspecting the proof of~\cite[Theorem 3.1]{Fang:1994wz}, one can observe that although the inequalities need only be satisfied locally on a small ball around a point $B(u_0,r)$, one must be able to choose the size of the ball and the constant $M$ in \ref{hol1} \emph{uniformly}, i.e. independently of the point $u_0$. The condition therefore requires at least a uniform structure, like that provided by a distance, to speak uniformly of the size of neighbourhoods. However, to our knowledge, there is no consensus on how to formulate a meaningful \emph{coordinate-invariant} version of~\ref{hol1}, and the authors do not define it in the statement of their theorem.
\medskip

A further shortcoming of the proposed extension of~\cite[Theorem 1.7, Theorem 3.1]{Fang:1994wz} to Riemannian manifolds that provides additional motivation for our investigation is the definition of the approximate Morse index given in~\cite[p.1599]{Fang:1994wz}. Their definition uses \emph{chart dependent} quantities and eliminates this dependence by taking the inf-sup over all charts. In addition to not exploiting in any way the assumed Riemannian structure, as stated, it appears unclear that such a quantity is accessible enough in any applications to be useful.

The underlying reason for such a definition is the fact that on a general manifold (without additional structure) there is no coordinate invariant analogue of the Hessian of a scalar function $\phi$ at arbitrary (non-critical) points. However, in the presence of a Riemannian structure, one \emph{can} define a covariant Hessian, simply as the covariant derivative $\nabla d \phi$ of the form $d\phi$ with respect to the Levi-Civita connection $\nabla$ of the Riemannan structure. More generally, a covariant Hessian can be defined whenever the manifold is equipped with a spray. We propose a new definition (see Definition~\ref{def: app morse}) of the approximate Morse index, based instead on this geometric structure. It is important to remark that the notion of Hessian is most useful when measuring the second order variation of $\phi$ along distinguished curves, the \emph{geodesics} of the spray.


\subsection*{Strategy of the proof of Theorem \ref{thm: application}: the use of two different sprays}

It follows from the above discussion that more work is required in order to generalise~\cite[Theorem 3.1]{Fang:1994wz} in a useful way to general Riemannian manifolds, if this is indeed possible at all. Our aim in the present paper is to make a first, limited, step in that direction, by developing an approach to the ideas in~\cite{Fang:1994wz} based on the notions of spray and geodesics. However, for reasons we shall explain below, we shall do this in the restricted setting introduced by Berestycki and Lions \cite[Section 8]{BerestyckiLions1983-2} where we study a functional on the constraint given by the submanifold $\M$ (see~\eqref{def-M}). This already provides an interesting framework that allows to handle several applications. Working in the presence of an ambient space also has the advantage of providing a univocal and convenient notion of Hölder continuity through the natural induced distance on $\M$.

The geometric structure at the heart of our considerations is that of a \emph{spray}. This general notion provides us with three important tools to explore the submanifold $\M$; a distinguished set of curves on $\M$ -- its geodesics --, a natural way to parallel transport vectors, and finally, a spray associates to any functional $\phi$ a natural second-order differential quantity we denote by $D^2\phi$, which is the covariant derivative of the one-form $\textrm{d}\phi$ restricted to $\M$. Just as the possibility to use information on the second-order derivative to perform deformations that decrease the value of the functional is key to the proof of \cite[Theorem 1.7 and Theorem 3.1]{Fang:1994wz}, we will do this with information on $D^2\phi$. This is accomplished in Lemma~\ref{lemm:3.3},  however it is more involved than the analogous statement in~\cite{Fang:1994wz}, as we are naturally confronted with the effects of curvature, for example, in the form of holonomy.

Although the general philosophy of our method is, \emph{in principle}, applicable to more general submanifolds than those we consider, it still has the disadvantage of requiring rather detailed knowledge about the geodesics; curves that are defined by a \emph{non}-linear second-order ordinary differential equation.  More precisely, in view of our result, the essential object we need to be able to understand in a \emph{quantitative} way is the exponential map \enquote{$\exp$} of the spray. In general, it is only defined, but nevertheless well understood, locally.  In particular, $\exp$ enables the construction of small convex neighbourhoods~\cite[Theorem 5.7, Chapter VIII \S 5]{Lang:1995ve} around each point and the \enquote{size} of the neighbourhood on which it is defined is related to the existence time of solutions. However, as we have mentioned, the nature of our result requires a certain uniformity (on a bounded set) that is slightly beyond these local considerations.  

A novel point that arises in our setting, illustrating a certain flexibility in our approach, is that we will work with two distinct, but nevertheless closely related, sprays. The first is chosen so that we have a nice manageable second-order quantity $D^2\phi$ from which it is effectively possible to extract interesting information, and the second being induced from the structure of $E$, we will develop this point further in Section~\ref{M-theorem}.
\smallskip

In the last part of the paper we present, in Theorem \ref{Objective1}, a version of Theorem \ref{thm: application} in which it is assumed that the family $\Phi_{\rho}: E \rightarrow \R$ consists of symmetric functionals. Namely: $\Phi_{\rho}(-x) = \Phi_{\rho}(x)$ for any $x \in E, \rho \in I$. Theorem \ref{Objective1} can be obtained as a direct consequence of Theorem \ref{M-theoremSym} which extends our key technical result, Theorem \ref{M-theorem}, stated in the next section. In Theorem \ref{M-theoremSym}, with the goal of deriving multiplicity results in the spirit of the symmetric mountain-pass theorem or of the fountain theorem,  we consider homotopic families which exhibit this symmetry. The main \emph{additional} difficulty is that we need to make deformations that preserve the symmetry of the homotopic family. 

\begin{theorem}\label{Objective1}
Let $I\subset (0,\infty)$ be an interval and consider a family of $C^2$ functionals $\Phi_\rho\colon E\to \R$ of the form 
$$\Phi_\rho(u)=A(u)-\rho B(u)$$
where $B(u)\geq 0$ for all $u\in E$ and 
$$
\text{either $A(u) \to +\infty$~ or $B(u) \to +\infty$ ~ as $u \in E$ and $\|u\| \to +\infty$.}$$

Suppose that $\Phi_\rho|_{S_\mu}$ is even for every $\rho \in I$, and moreover that $\Phi'_\rho$ and $\Phi''_\rho$ are $\alpha$-H\"older on bounded sets in the sense of Definition~\ref{Holder continuous}~(\ref{hol2}) for some $\alpha \in (0,1]$.
Finally, suppose that there exists a $N\in \N$, $N \geq 2$, given odd functions $\gamma_i : S^{N-2} \mapsto S_{\mu}$ where $i =0,1$, such that the set 
\begin{equation}\label{SMP1}
\Gamma_N=\{ \gamma\in C([0,1]\times S^{N-2}, S_\mu)\colon \gamma(t,\cdot) \ \text{is odd }, \gamma (0, \cdot) = \gamma_0, \gamma(1, \cdot) = \gamma_1 \}
\end{equation}
is non-empty and that 
\begin{equation}\label{SMP2}
c_{\rho,N}=\inf_{\gamma\in\Gamma}\max_{[0,1] \times S^{N-2}}\Phi_\rho(\gamma(t,x))> \max_{x \in S^{N-2}}\big\{\Phi_\rho (\gamma_0(x)), \Phi_\rho (\gamma_1(x))\big\}, \quad \forall \rho \in I.
\end{equation}
Then, for almost every $\rho \in I$, there exist sequences $\{u_n\} \subset \M$ and $\zeta_n \to 0^+$ such that, as $n \to + \infty$,
\begin{itemize}
\item[(i)] $\Phi_\rho(u_n) \to c_{\rho,N}$;
\item[(ii)] $||\Phi'_\rho|_{\M}(u_n)||_* \to 0$;
\item[(iii)] $\{u_n\}$ is bounded in $E$;
\item[(iv)] $\tilde m_{\zeta_n}(u_n) \leq N$.
\end{itemize} 
\end{theorem}

\begin{remark}
The bound $\tilde m_{\zeta_n}(u_n) \leq N$, where one could expect $\tilde m_{\zeta_n}(u_n) \leq N-1$, arises from the fact that $[0,1]\times S^{N-2}$ need to be viewed as a subset of $\R^N$.
In the case $N=2$, we can improve (iv) to $\tilde m_{\zeta_n}(u_n) \leq N-1$ through the identification of $[0,1]\times \{-1,1\}$ with a subset of $\R$.  However, note that (iv) provides only an inequality and that, in most applications of Theorem \ref{Objective1}, what shall be used is that the approximate Morse index of $\{u_n\} \subset \M$ is uniformly bounded.
\end{remark}

\begin{remark} Theorem \ref{Objective1} is only one of the possible corollaries of Theorem \ref{M-theoremSym}. Its interest lies in the fact that, in some applications, one may hope to check that, for several values of $N \in \N$, the class $\Gamma_N$ is not empty and that~\eqref{SMP2} holds. Assuming in addition that  $c_{\rho, N}\neq c_{\rho,M}$ for at least one pair $(N,M)\in \N^2$, this would lead to the existence of Palais-Smale sequences satisfying (i)-(iv) at distinct levels $c_{\rho,N}$ and open the way to establishing the existence of multiple critical points. Actually, this type of strategy was exploited in \cite{BaVa} where, on the specific problem they consider there, the authors obtain the existence of infinitely solutions with a prescribed $L^2$ norm. This paper was one of our motivations for formulating Theorem \ref{Objective1} which we hope will lead, in particular, to multiplicity results on graphs, extending the existence results of \cite{BoChJeSo:2022,ChJeSo:2022}. 
\end{remark}

\medskip 

The paper is organised as follows. In Section \ref{Theorem-FGrevisited} we present our key technical result, Theorem \ref{M-theorem} below, which is based on the deformations constructed in Lemma \ref{lemm:3.3}. Theorem \ref{thm: application} is proved in Section \ref{Theorem-application} where we also explore some of its consequences when the sequence $\{u_n\} \subset \M$ converges. In particular, we derive information on the Morse index of the limiting critical point, see Theorem \ref{thm: limit}. Section \ref{Theorem-Symmetric} is devoted to establishing Theorem \ref{M-theoremSym} which extends Theorem \ref{M-theorem} to a symmetric setting. 

\subsection*{Acknowledgments and remark on the numbering} 
We thank the anonymous referee for carefully reading our original manuscript, and for their suggestions for possible extensions, which have lead to this richer, extended version. As a result, we had to slightly reconsider the organisation of the paper. In particular, the numbering of the mathematical environments (theorems, lemmata, \ldots) has changed with respect to the first version of the paper. Since some results have been already quoted in the literature (e.g. in \cite{BoChJeSo:2022,ChJeSo:2022}), we prefer to explicitly mention that, for instance:
\begin{itemize}
\item Theorem \ref{thm: 0-application} corresponds to Theorem 1 of the first version; 
\item Theorems \ref{thm: application} and \ref{Objective1} are new; 
\item Theorem \ref{M-theorem} corresponds to Theorem 2 in the first version; 
\item Remarks \ref{Gradient}-\ref{space-of-functions} correspond to Remark 1.3-1.5 in the first version; 
\item Corollary \ref{coro:step1} corresponds to Corollary 1 in the first version.
\end{itemize}

\section{Statement and proof of Theorem \ref{M-theorem} }\label{Theorem-FGrevisited}

This section, which is the heart of the paper, is devoted to the proof of Theorem \ref{M-theorem} below.
\newline
We will work quite intensively with the induced topology on $\M$. In order to simplify notations, we will write $B(\M; u_0 ,\delta)=B(u_0,\delta)\cap \M$ for any $u_0\in \M$.



\begin{theorem}\label{M-theorem}
Let $\phi$ be a $C^2$-functional on $E$, satisfying \ref{hol2} for some $\alpha \in (0,1]$; let $\mathcal{F}$ be a homotopy-stable family of $\M$ of dimension $n$ with boundary $B$, such that $$c:= \inf_{A \in \mathcal{F}} \max_{u \in A}\phi(u) > \max_B \phi$$ is finite.
Let $R>1$, $0< \alpha_1 \leq \frac{\alpha}{2(\alpha +2)}<1$ and $\va>0$.
Then for any $A \in \mathcal{F}$ with $\max_{u \in A} \phi(u) \leq c + \varepsilon$ satisfying 
\begin{equation}\label{top-bounded}
K:= \{ u \in A \,|\, \phi(u) \geq c - \va \} \subset B(0,R-1), 
\end{equation}
assuming that $\varepsilon >0$ is sufficiently small there exists $u_{\varepsilon} \in \M$ such that
\begin{enumerate}
\item $c- \va \leq \phi(u_{\va}) \leq c + \va$ ;
\item $||\phi'|_{\M}(u_{\va})|| \leq 3 \va^{\alpha_1}$ ;
\item $ u_{\va} \in A$ ;
\item If $D^2\phi(u_{\va}) [w, w]  < -  \va^{\alpha_1} \|w \|^2$ for all $w \neq 0$ in a subspace $W$ of $T_{u_{\varepsilon}}\M$, then $\dim \, W \leq n$.
\end{enumerate}
\end{theorem}

In order to prove Theorem \ref{M-theorem} we recall, for the convenience of the reader, some results presented in \cite{Fang:1994wz}; the first is \cite[Lemma 2.5]{Fang:1994wz}.

\begin{lemma}\label{Lemma 2.5-FG}
For each $n \in \mathbb{N}^*,$ there is an integer $N(n) \leq (2 \sqrt{n+1}+ 2)^n$ such that for any compact subset $C \subset \R^n$ and any $\va >0$, there exist a finite number of distinct points $\{x_i : 1 \leq i \leq k \}$ with the following properties :
\begin{enumerate}
\item $ \displaystyle C \subset \bigcup_{i=1}^k B(x_i, \frac{\va}{4})  \subset N_{\frac{\va}{2}}(C)$;
\item The intersection of any distinct $N(n)$ elements of the cover $( \overline{B}(x_i, \frac{\va}{2}))_{i=1}^k$ is empty.
\end{enumerate}
\end{lemma}


The following result is \cite[Corollary 2.4]{Fang:1994wz}.

\begin{lemma}\label{Corollary 2.4-FG}
Let $V$ be a Banach space, let $K$ be a closed subset of $\mathbb{R}^n$, and let $\theta$ be a continuous mapping from $K$ into the unit sphere $S_V$ of $V$. If $n < \dim V$, then $\theta$ can be extended to a continuous mapping from $\mathbb{R}^n$ into $S_V.$ 
\end{lemma}

\subsection{The geodesic spray coming from the scalar product of $H$}
We will now introduce the geometric structure at the heart of our discussion: sprays.
These are second-order vector fields $F$ -- a specific type of vector field on the tangent bundle $TX$ of a manifold $X$ -- that are locally determined by a smooth family of quadratic forms on the model space $E$. Its geodesics are then defined as the projections of the integral curves of $F$ to $X$. 

In our setting, we shall consider sprays on $\M$ that are restrictions of sprays on $E\setminus\{0\}$, which can be completely described by a map: $F : T(E\setminus\{0\})\cong E\setminus\{0\}\times E \rightarrow TT(E\setminus\{0\}) \cong  E\setminus\{0\}\times E \times E\times E$ of the form:
\[ F(x,v)=(x,v,v,f(x,v)), \]
where $f(x,\cdot)$ is a quadratic form on $E$.
Geodesics are then projections onto the first factor of solutions $(\gamma, v)$, to the system:
\[(\gamma,v,\dot{\gamma},\dot{v})=F(\gamma,v),\]
which of course amounts to the second order ODE:
\[\ddot{\gamma}=f(\gamma,\dot{\gamma}). \] 
We refer the reader to \cite{Lang:1995ve} for a detailed general discussion.

Whilst $\M=\{ u \in E, (u,u)=|u|^2=\mu\}$ naturally inherits from $E$ the structure a smooth Riemannian manifold, the scalar product $(\cdot, \cdot)$ of $H$, on the other hand, does not necessarily equip it with a Riemannian structure: the tangent spaces, in general, will not be complete under the induced norm. Nevertheless, it does equip $\M$ with a spray.
This can be described as follows, let $F_1: E\setminus\{0\} \times E \rightarrow E\setminus\{0\}\times E \times E \times E$ be defined by:
\begin{equation}\label{SprayF1} F_1(u,v) = (u,v,v,-\frac{|v|^2}{|u|^2}u),\end{equation}
it is clear that this is a spray on $E\setminus\{0\}$, furthermore:
\begin{lemma}
\label{lemm:restrict_spray}
The spray $F_1$ restricts to a spray on $\M$, i.e. $F_1: T\M \rightarrow TT\M$.
\end{lemma}
\begin{proof}
This is a direct consequence of the following characterisations of $T\M$ and $TT\M$:
\[ \begin{gathered} T\M=\{ (u,v) \in TE, u\in \M, v \in T_u\M\}=\{ (u,v) \in TE, u\in \M, (u,v)=0\}, \\ TT\M = \{(u,v,v',w)\in TTE, u\in \M, (u,v)=0, (u,v')=0, (u,w)+(v',v)=0\}.  \end{gathered} \]

These characterisations result directly from the fact that $T\M$ is defined in $TE=E\times E$ by the equations:
\[ (u,v) \in T\M \Leftrightarrow u \in \M, v \in T_u\M \Leftrightarrow (u,v)\in E\times E, |u|^2=\mu, (u,v)=0.\] Consequently, if $\psi : E\times E \rightarrow \mathbb{R}^2$ is defined by $\psi(u,v)=(|u|^2 - \mu,(u,v))$, $\psi^{-1}(\{(0,0)\})=T\M$ and hence $T_{(u,v)}T\M= \ker \psi'(u,v).$
\end{proof}

The geodesics of $F_1$ satisfy the equation:
\begin{equation}\label{curve_eq2} \sigma''(t)+ \frac{|\sigma'(t)|^2}{|\sigma(t)|^2}\sigma(t)=0. \end{equation} 

\begin{remark} Lemma~\ref{lemm:restrict_spray} is the differential expression of the fact that if the initial data of Equation~\eqref{curve_eq2} is taken in $T\M$ then the solution lives on $\M$, i.e. for any $t\in I$, $|\sigma(t)|^2=\mu$.\end{remark} The principal particularity of this spray is that: 
\begin{lemma}\label{sprayF1conservesH}Let $\sigma: I \rightarrow \M$ be a geodesic of the spray $F_1$ on $\M$, then: \[ \forall t\in I,\quad |\sigma'(t)|^2=|\sigma'(0)|^2.\] \end{lemma}
A major advantage of the spray $F_1$ is that the geodesics, with initial data $(u,v)\in T\M$, are explicitly known:
\begin{equation} \label{eq:exp_explicite} \sigma(t,u,v)=\cos(\omega t)u + \frac{\sin(\omega t)}{\omega}v, \qquad \omega=\frac{|v|}{\sqrt{\mu}},\quad t\in \mathbb{R}.\end{equation}
These curves are globally defined and, as we observed in Lemma~\ref{sprayF1conservesH}, have constant speed with respect to the norm of $H$, but \emph{not} with that of $E$. They also possess the following homogeneity property of geodesics $a\in \mathbb{R}$:
\[ \sigma(at,u,v)=\sigma(t,u,av). \]

Let us introduce and state a few results about the exponential map associated with the spray:
\begin{lemma} \label{lemm:exp}
For $(u,v)\in T\M$, define: \[\exp_u(v)=\sigma(1,u,v),\] and $\exp: T\M \rightarrow \M \times \M$, \[\exp(u,v)=(u,\exp_u (v)).\]
Then:
\[ \textrm{d}\exp_{(u,v)}(v',w)=\underbrace{\left( v', \cos\left(\frac{|v|}{\sqrt{\mu}}\right)\left[v' +\frac{(w,v)}{|v|^2}v\right]+\sin\left(\frac{|v|}{\sqrt{\mu}}\right)\sqrt{\mu}\left[ \frac{w}{|v|} - \frac{(w,v)}{|v|^3}v - \frac{(w,v)}{\mu |v|} u\right]  \right)}_{ \in T_{(u,\exp_u(v))}\M\times\M\cong T_u\M\times T_{\exp_u(v)}\M.}\]
Let $\mathcal{U}=\{(u,v)\in T\M, |v| < \sqrt{\mu}\pi\}$, then $\exp \mathcal{U}$ is open and $\exp$ is a diffeomorphism onto its image, which is $\{(u_1,u_2)\in \M\times \M, u_2 \neq -u_1\}$.
For fixed $u_0\in \M$: \begin{equation} \label{eq:expinv} \exp_{u_0}^{-1}(u) =\frac{\arccos\frac{(u,u_0)}{\mu}}{\sqrt{1-\frac{(u,u_0)^2}{\mu^2}}}\left(u-\frac{(u,u_0)}{\mu}u_0\right), \quad u\neq -u_0 \end{equation} and:
\begin{equation} \label{eq:l2normexpinv} |\exp_{u_0}^{-1}(u)|= \sqrt{\mu} \arccos \frac{(u,u_0)}{\mu}.  \end{equation} 

\end{lemma}
\begin{proof}
$\mathcal{U}=\{(u,v)\in TE=E\times E, |v|< \sqrt{\mu\pi} \}\cap T\M$ is open in $T\M$ and, given the expression of the differential, one can see that it is (continuously) invertible at each point in $\mathcal{U}$. Let us briefly justify Equation~\eqref{eq:expinv}, we solve directly for $v$:
\[u=\cos(\frac{|v|}{\sqrt{\mu}})u_0 + \sinc(\frac{|v|}{\sqrt{\mu}})v, \qquad v\in T_{u_0}\M, |v| < \sqrt{\mu}\pi, \]
where $u\neq -u_0$. Since $v\in T_{u_0}\M$, it follows that:
\[(u,u_0)=\mu\cos(\frac{|v|}{\sqrt{\mu}}) \Rightarrow |v| = \sqrt{\mu}\arccos \frac{(u,u_0)}{\mu}.\]
Then:
\[\sinc \frac{|v|}{\sqrt{\mu}}= \left\{ \begin{array}{lc}\sqrt{1-\frac{(u,u_0)^2}{\mu^2}}\left(\arccos \frac{(u,u_0)}{\mu}\right)^{-1} & \textrm{if $u\neq u_0,$} \\ 1 & \textrm{if $u =u_0$}.\end{array} \right.\]
Reinjecting these expressions into the original equation provides the expression given by~\eqref{eq:expinv} for $v$.
\end{proof}

\begin{remark}\label{homogeneity}
It will be convenient to note that, due to the homogeneity property, for any $t\in \mathbb{R}$, \begin{equation} \label{rema:exp_sigma}\exp_u(tv)=\sigma(1,u,tv)=\sigma(t,u,v). \end{equation}
\end{remark}
\subsection{Geodesics and $D^2\phi$}\label{relevanceD2}
We shall now explain the relevance of Definition~\ref{def D}. In short, it is the covariant derivative of $\textrm{d}\phi$ (viewed as a one-form) on $\M$ induced by $F_1$, but we shall not develop this point of view further. For us, the important point is given by the following lemma.
\begin{lemma}\label{LinkHessianGeodesic}
Let $\phi$ be a $C^2$-functional on $E$ and $\sigma: (-\varepsilon,\varepsilon) \rightarrow E$ an arbitrary curve satisfying Equation~\eqref{curve_eq2}, then for any $t\in (-\varepsilon,\varepsilon)$:
\[ \frac{\textrm{d}^2}{\textrm{d}t^2}\phi(\sigma(t))=D^2\phi(\sigma(t))[\sigma'(t), \sigma'(t)]. \]
\end{lemma}
\begin{proof}
\[\frac{\textrm{d}}{\textrm{d}t}\phi(\sigma(t))=\phi'(\sigma(t))\cdot\sigma'(t). \]
Therefore:
\[\begin{aligned} \frac{\textrm{d}^2}{\textrm{d}t^2}\phi(\sigma(t))&= \phi''(\sigma(t))[\sigma'(t),\sigma'(t)]+ \phi'(\sigma(t))\cdot\sigma''(t)\\&= \phi''(\sigma(t))[\sigma'(t),\sigma'(t)]-\frac{|\sigma'(t)|^2}{|\sigma(t)|^2}\phi'(\sigma(t))\cdot \sigma(t). \end{aligned}\]
\end{proof}

\begin{remark}\label{rema:hessian}
If $u$ is a critical point of the functional $\phi|_{\M}$, then the above also shows that the restriction of $D^2\phi(u)$ to $T_u\M$ coincides with the \emph{Hessian} of $\phi|_{\M}$ at $u$ (as defined, for instance, in~\cite[p307]{PALAIS1963299}); this follows, for instance, from~\cite[Proposition 5.2.3]{Ma}. This is also a direct consequence of the fact that the restriction of $D^2\phi(u)$ to $T_u\M$ is the covariant derivative associated with the spray $F_1$ of the one-form $d\phi|_{\M}$.
\end{remark}

\begin{lemma}
\label{lemm:holder}
Let $\phi$ be a $C^2$-functional on $E$, $\alpha \in (0,1]$ such that $\phi, \phi''$ are $\alpha$-Hölder continuous on bounded sets (\ref{hol2}). 
Then for any $R\geq 0$, one can find $C\in \mathbb{R}_+^*$ such that for any $u_1,u_2 \in \M\cap B(0,R)$: 
\[\begin{cases} ||D^2\phi(u_1)-D^2\phi(u_2)||_{**}\leq C ||u_1-u_2||^\alpha & \textrm{if $||u_1-u_2||\leq 1$}, \\
||D^2\phi(u_1)-D^2\phi(u_2)||_{**}\leq C ||u_1-u_2|| & \textrm{if $||u_1-u_2||> 1$}. \end{cases} \]
 \end{lemma} 
\begin{proof}
Let us estimate $|\phi'(u_1)u_1-\phi'(u_2)u_2|$:
\[
\begin{aligned}
|\phi'(u_1)u_1-\phi'(u_2)u_2| &\leq |(\phi'(u_1)-\phi'(u_2))\cdot u_1| + |\phi'(u_2)\cdot(u_2-u_1)| \\&\leq M R ||u_1-u_2||^\alpha + MR^{\alpha}||u_1-u_2||.
\end{aligned}
\]
It then follows that:
\[||D^2\phi(u_1)-D^2\phi(u_2)||_{**} \leq M\left(1+ \frac{R}{\mu}\right)||u_1-u_2||^\alpha+ \frac{MR^{\alpha}}{\mu}||u_1-u_2||. \qedhere\]
\end{proof}

\begin{corollary}
\label{coro:step1}
Let $\phi$ be a $C^2$-functional on $E$ such that $\phi'$ and $\phi''$ are $\alpha$-Hölder continuous on bounded sets. Assume that for some $R>1$, $u_0 \in \M \cap B(0,R-1)$, one can find a constant $\beta \in (0,1)$ and a subspace $W\in T_{u_0}\M$ such that for any $w\in W:$
\[ D^2\phi(u_0)[w,w] < -\beta ||w||^2. \]
Then, for any $u\in B(\M; u_0,\delta_1)$ and $w \in W$:
\begin{equation}\label{conditionD2step1} D^2\phi(u)[w,w] < - \frac{3\beta}{4} ||w||^2,\end{equation}
where
\begin{equation}\label{def:delta1}
\delta_1 = \min \left(1,\left(\frac{\beta}{4C}\right)^{\frac{1}{\alpha}}\right), \quad \mbox{with} \quad C= M\left(1+ \frac{R+R^{\alpha}}{\mu}\right).
\end{equation}
\end{corollary}
\begin{proof}
This follows directly from Lemma~\ref{lemm:holder} since if $u\in B(u_0,\delta),$ $\delta\leq 1$, $u\in B(0,R)$. 
\end{proof}

\noindent
For future reference, let us fix 
$K(R) \geq 1$ such that
\begin{equation}\label{def:K}
||D^2 \phi (u) ||_{**} \leq K(R) \quad \mbox{and} \quad ||\phi'(u)||_*\leq K(R) \quad \mbox{for any } u \in B(0,R)\cap \M.
\end{equation}

\subsection{Transporting $W$ to other points}
The condition expressed in Equation~\eqref{conditionD2step1} in the result of Corollary~\ref{coro:step1} is not an intrinsic condition on $\M$. The problem is that in general, $W$ is not tangent to the manifold at points $u\neq u_0$ on $\M$. To overcome this, $W$ should also move with the point to remain tangent to the manifold. In other words, for each point $u$ in a small neighbourhood of $u_0$, we should exhibit a subspace $W(u)$ of $T_uS_{\mu}$, which satisfies a condition of the same nature as~\eqref{conditionD2step1}.

Our approach to this will be based on the notion of \emph{parallel transport} of vectors along curves associated, however, to a \emph{second} spray $F_2\neq F_1$ . Parallel transport of a given vector $w_0\in T_{u_0}\M$ along geodesics of $F_1$ originating from $u_0$, we will obtain a vector field on a ball $B(\M; u_0,\delta)$. 
The spray $F_2$, introduced below, is, contrary to $F_1$, naturally associated with the Riemannian structure on $\M$ induced by $E$; it is known as the \emph{canonical spray}, we refer the reader to~\cite[Chapter VII, \S 7, Chapter VIII, \S 4]{Lang:1995ve}) for the detailed construction.

To describe the spray $F_2$ in our setting, introduce the injective linear map $G: E\rightarrow E$,  such that for any $u\in E$, $Gu$ is the unique vector in $E$ that satisfies:
\begin{equation}\label{definition-G}
\forall h \in E, (u, h) = \langle Gu,h\rangle. 
\end{equation}
Note that:
\[ ||Gu||^2= \langle Gu, Gu \rangle = (u,Gu) \leq ||Gu|| \cdot |u|, \]
therefore, for $u\neq 0$,
\[ ||Gu|| \leq |u| \leq ||u||. \]
Furthermore, if $u\in \M$ then $\langle Gu, u \rangle = |u|^2=\mu$, thus $\mu \leq ||Gu||\cdot ||u||$. We collect these observations in the following lemma for future reference:
\begin{lemma}
The injective linear map $G$ has norm $||G|| \leq 1$ and for any $u\in\M:$
 \begin{equation}\label{eq:lbG} ||Gu|| \geq \frac{\mu}{||u||}.\end{equation}
\end{lemma} 
In terms of $G$, the canonical spray $F_2 : T(E\setminus\{0\}) \rightarrow TT(E\setminus\{0\})$ on $E\setminus\{0\}$ can be defined by: \begin{equation}\label{SprayF2} F_2(u,v)=f(u,v,v,-|v|^2\frac{Gu}{||Gu||^2}).\end{equation}
One can check in a straightforward manner that it restricts to a spray on $\M$; its particularity is that it preserves the induced Riemannian structure on $\M$ (see Lemma~\ref{para tran} below.) 

On this specific example, let us describe how a spray is used to define parallel transport of vectors along curves. The discussion applies to arbitrary curves, however in our situation, they will systematically be geodesics of the spray $F_1$, more precisely: let $u_0 \in \M$ and $u \in B(\M; u_0,\delta)$; assume $\delta <2\sqrt{\mu}$ so that $-u_0\notin B(\M;u_0,\delta)$. Set $v=\exp_{u_0}^{-1}(u)$ and\footnote{See \eqref{rema:exp_sigma} in Remark~\ref{homogeneity}.}, for any $t\in[0,1]$, $\sigma(t)=\exp_{u_0}(tv)=\sigma(t,u_0,v)$. If $w_0 \in T_{u_0}S$, we define a vector field on $S_\mu$ along the curve $\sigma$ by the following differential equation constructed from the symmetric bilinear operator associated with the spray $F_2$:
\begin{equation}\label{eq://tr} \varphi'(t) + (\varphi(t), \sigma'(t))\frac{G\sigma(t)}{||G\sigma(t)||^2} =0, \quad \varphi(0)=w_0. 
\end{equation}

We quote the following properties of parallel transport:
\begin{lemma}\label{para tran}
Let $u_0, u \in \M$ such that $v=\exp_{u_0}^{-1}(u)$ and define $\sigma(t)=\exp_{u_0}(tv), t\in [0,1]$. Define a map: $T_{u_0,u,\sigma}: T_{u_0}\M \rightarrow T_{u}\M$ by $T_{u_0,u,\sigma}w_0= \varphi(1)$ where $t\mapsto \varphi(t)$ is the solution to Equation~\eqref{eq://tr}.
Then:\begin{itemize} \item $T_{u_0,u,\sigma}$ is a linear isometry (for the induced Hilbert space structure on the tangent spaces), 
\item The inverse is given by: $T^{-1}_{u_0,u,\sigma}= T_{u,u_0,\tilde{\sigma}}$ where $\tilde{\sigma}(t)=\sigma(1-t), t\in [0,1]$.
\end{itemize}
\end{lemma}
\begin{proof}
We shall only prove that it is an isometry; differentiating $(\varphi,\sigma)$:
\[ (\varphi'(t),\sigma(t)) + (\varphi(t),\sigma'(t))= -\frac{(G\sigma(t),\sigma(t))}{||G\sigma(t)||^2}(\varphi(t),\sigma'(t)) + (\varphi(t),\sigma'(t))=0.\]
Since $(G\sigma(t),\sigma(t))=||G\sigma(t)||^2.$
Similarly differentiating $\langle \varphi(t), \varphi(t)\rangle$ yields:
\[
\begin{split} 
\langle \varphi(t), \varphi(t)\rangle'= 2\langle \varphi'(t), \varphi(t)\rangle &= -2\frac{(\varphi(t),\sigma'(t))}{||G\sigma(t)||^2}\langle G\sigma(t),\varphi(t) \rangle \\
& = -2\frac{(\varphi(t),\sigma'(t))}{||G\sigma(t)||^2}( \sigma(t),\varphi(t))=0.  \end{split} 
\]
\end{proof}

Since all our curves will be geodesics of $F_1$ that originate from a fixed point $u_0\in \M$, we abbreviate: $T_{u_0,u,\sigma}$ to $T_u$ and say that we \emph{parallel transport radially} from $u_0$. Given $w_0 \in W$, we shall define a vector field on $B(\M; u_0,\delta)$ by:
\[w(u)=T_u w_0, \quad u \in B(\M; u_0, \delta).\]
The spray allows us to compare vectors at different points of $\M$ and, the vector field defined by Equation~\eqref{eq://tr} is, relative to the spray $F_2$, \enquote{constant} along the curve $\sigma$. However, viewing things in $E$, it must change to remain tangent to $\M$; the next lemma estimates this change. 
\begin{lemma}
\label{lemm:diff_transport}
Let $w_0 \in T_{u_0} \M$, $u_0\in B(0,R-1) \cap \M$, $u \in B(\M; u_0,\delta_0)$ with $\delta_0 :=\min(1, \sqrt{\mu})$. Then, if $v=\exp_{u_0}^{-1}(u)$: \[ || T_u w_0 - w_0|| \leq C(R, \mu) \, |v|\cdot ||w_0||  \]
where we have set
\begin{equation}\label{def:C(R,mu)}
C(R, \mu) = \frac{1}{\mu} \big(R+\frac{(R-1)}{\sqrt{\mu}}\big).
\end{equation}
\end{lemma}
\begin{proof}
By definition of $T_u w_0$, it is sufficient to estimate $||\varphi(1)-\varphi(0)||$, where $\varphi$ is defined by Equation~\eqref{eq://tr}. By the mean-value theorem one has:
\[\begin{aligned} ||\varphi'(t)|| \leq \frac{||\varphi(0)|| \cdot |\sigma'(0)|}{||G\sigma(t)||} \leq \frac{||w_0|| \cdot |v| \cdot ||\sigma(t)||}{\mu}, \end{aligned} \]
where we have used Equation~\eqref{eq:lbG} and the fact that $|\sigma'(t)|=|\sigma'(0)|=|v|$ by Lemma~\ref{sprayF1conservesH}.

Now, taking into account that $\delta_0 \leq \sqrt{\mu}$ guarantees that $(u, u_0) \geq 0$, we estimate $||	\sigma(t)||$ as follows:
\[\begin{aligned} ||\sigma(t)|| &\leq R-1 + \sqrt{\mu} \sin\left(\frac{|v|}{\sqrt{\mu}}t\right) \frac{||v||}{|v|} \\ &\leq R-1+\sqrt{\mu}\sin\left(\frac{|v|}{\sqrt{\mu}}\right) \frac{||v||}{|v|}. \end{aligned}  \]
From Equations~\eqref{eq:expinv}, \eqref{eq:l2normexpinv} from Lemma~\ref{lemm:exp}: \[ |v|=\sqrt{\mu}\arccos\frac{(u,u_0)}{\mu}, \frac{v}{|v|}=\frac{\sqrt{\mu}}{\sqrt{\mu^2-(u,u_0)^2}}(u-\frac{(u,u_0)}{\mu}u_0). \]
Therefore: 
\[\begin{aligned} ||\sigma(t)|| &\leq R-1 + ||u-\frac{(u,u_0)}{\mu}u_0||\\&= R-1+||u-u_0 -\frac{(u-u_0,u_0)}{\mu}u_0|| \\&\leq R-1 + \delta_0+ \frac{(R-1)\delta}{\sqrt{\mu}}.  \end{aligned} \]
The result follows as $\delta_0 \leq 1$.
\end{proof}

We can now prove the main result of this section:
\begin{lemma}
\label{lemm:step1}
Let $\phi: E \rightarrow \mathbb{R}$ such that $\phi'$ and $\phi''$ are $\alpha$-Hölder continuous on bounded sets, $u_0\in B(0,R-1) \cap \M$, $R>1$, and suppose that one can find $\beta\in (0,1)$ and a subspace $W\equiv W(u_0)\subset T_{u_0}\M$ of dimension $n+1$ such that:
\[D^2\phi(u_0)[w,w] < -\beta ||w||^2.\]
For any $u\in B(\M; u_0, 2\sqrt{\mu})$ set $W(u)=T_uW$ the image of $W$ under radial parallel transport, then for any $u\in B(\M; u_0, \delta_2)$ and $w \in W(u)$,
\[D^2\phi(u)[w,w] < -\frac{\beta}{2} ||w||^2, \]
where 
\begin{equation}\label{def:delta2}
\delta_2 = \begin{cases}  \displaystyle \min\left(\sqrt{\mu}\left(1-\cos\left(\frac{\beta}{8 \sqrt{\mu} \,K(R) \ C(R, \mu)}\right)\right),\delta_1, \delta_0\right) & \textrm{if  } \displaystyle \frac{\beta}{8 \sqrt{\mu} \,K(R) \ C(R, \mu)}<\pi,\\ \min(\delta_1, \delta_0) & \textrm{otherwise}\end{cases}
\end{equation}
with $\delta_0 >0$ being defined in Lemma~\ref{lemm:diff_transport}, $\delta_1 >0$ in Corollary~\ref{coro:step1}, $K(R)>0$ and $C(R, \mu) >0$ in \eqref{def:K} and  \eqref{def:C(R,mu)} respectively.
\end{lemma}
\begin{proof}
According to Corollary~\ref{coro:step1}, one can already find $\delta_1 >0$ satisfying Equation \eqref{def:delta1} such that for any $u\in B(\M; u_0,\delta_1)$ and any $w\in W(u_0)$:
\[D^2\phi(u)[w,w] < -\frac{3\beta}{4}||w||^2.\]
Now,  fix $u\in B(\M; u_0,\delta_2)$, where $\delta_2 >0$ satisfies \eqref{def:delta2}, and choose $(\tilde{w}_i, \dots, \tilde{w}_n)$ an orthonormal basis of $T_uW$ that diagonalises the bilinear form $D^2\phi(u)$. Set for each $i\in \{1,\dots,n+1\}$, $w_i=T_u^{-1}\tilde{w}_i$; this defines an orthonormal basis of $W$. Using Lemma~\ref{lemm:diff_transport}, and the definition of $K(R)$, we have for any $i\in\{1,\dots,n+1\}$:
\[\begin{aligned} |D^2\phi(u)[w_i,w_i]- D^2\phi(u)[T_uw_i,T_uw_i]|&=|D^2\phi(u)[w_i-T_uw_i,w_i+T_uw_i]|\\ &\leq \left(R+\frac{R-1}{\sqrt{\mu}} \right) \frac{2K(R)}{\mu}|\exp^{-1}_{u_0}(u)|\\
&=2K(R)C(R,\mu)|\exp^{-1}_{u_0}(u)|.
\end{aligned}\]
Note that: \[\frac{(u,u_0)}{\mu} = \frac{(u-u_0,u_0)}{\mu}+1\geq 1-\frac{\delta_2}{\sqrt{\mu}}, \] thus, by Equation~\eqref{eq:l2normexpinv}: \[|\exp_{u_0}^{-1}(u)|\leq \sqrt{\mu}\arccos(1-\frac{\delta_2}{\sqrt{\mu}}).\]
Thus, for any $u\in B(\M; u_0,\delta_2)$, and any $i\in\{1,\dots,n+1\}$, 
 in view of the condition on $\delta_2 >0$ given in \eqref{def:delta2}, if $\frac{\beta}{8 \sqrt{\mu} \,K(R) \ C(R, \mu)}<\pi$, then 
\[
2K(R)C(R,\mu)|\exp^{-1}_{u_0}(u)|\le \frac{\beta}{4};
\] otherwise, 
\[
2K(R)C(R,\mu)|\exp^{-1}_{u_0}(u)|\le 2K(R)C(R,\mu)\frac{\sqrt{\mu}\pi}{2}\le \frac{\beta}{4}.
\]
Hence
\[ D^2\phi(u)[T_uw_i,T_uw_i] < -\frac{\beta}{2} = -\frac{\beta}{2}||T_uw_i||^2.  \]
Since $T_uw_i= \tilde{w}_i$ by definition, the inequality extends to all of $W(u)$ by bilinearity.
\end{proof}

\subsection{Fang-Ghoussoub's result}
The remainder of this section is devoted to the proof of Theorem~\ref{M-theorem}, which will be divided into several lemmata. Our initial goal is to show how to use the second-order information on a functional $\phi$ in the conclusion of Lemma~\ref{lemm:step1} to find appropriate deformations of an element $A$ of the $\leq d$-dimensional family $\mathcal{F}$ so as to locally decrease the value of the functional $\phi$.  As in~\cite{Fang:1994wz}, we shall probe $\M$ by running out along specific curves from points in a neighbourhood of a point $u_0\in A$; in our case these curves will be geodesics of the spray $F_1$.

Our first lemma identifies, uniformly on a bounded set, a maximal travel time that does not take us too far on $\M$.
\begin{lemma}
\label{lemm:tmax}
Let $\delta >0$ and assume $u_0 \in B(0,R)\cap \M $, $u\in B(\M; u_0,\frac{\delta}{2})$, $v\in T_{u}\M\setminus\{0\}$, $||v|| \leq 1$. 
Set $t_{max}=\min \{\frac{2\mu}{R}, \frac{\delta}{4} \}$, then, for any $0\leq t < t_{max}$:  \[\sigma(t,u,v)\in B(u_0,\delta).\]
\end{lemma}
\begin{proof}
Note that since $0<|v| \leq ||v|| \leq 1$ it follows that $\sqrt{\mu}\leq \frac{1}{\omega}$, so $ t_{max} \leq \frac{2\mu}{R} \leq \frac{2}{R\omega^2}$. Then:
\[\begin{aligned} ||\sigma(t,u,v)-u_0|| &< |\cos(\omega t)|\frac{\delta}{2} + |1-\cos(\omega t)|R + \frac{|\sin(\omega t)|}{\omega} \\ &< \frac{\delta}{2} + 2|\sin\frac{\omega t}{2}|\left( |\sin\frac{\omega t}{2}| R + \frac{1}{\omega} \right)\\ &< \frac{\delta}{2} +\frac{4}{\omega}|\sin \frac{\omega t}{2}| \leq \delta.\end{aligned}\]
\end{proof}

We will encounter two geometric issues in the sequel; the first comes from the fact that we use two different sprays. Since the geodesics of $F_1$ differ from those of $F_2$ they do not parallel transport their initial velocity vector according to the spray $F_2$; therefore in general even if $\sigma'(0,u,w)\in W(u)$, $\sigma'(1,u,w)\notin W(\sigma(1,u,w))$. The second is a global expression of curvature: parallel transport along a closed curve $\gamma: [0,1] \rightarrow \M$, $\gamma(0)=\gamma(1)=u_0$, defines a linear operator: $A: T_{u_0}\M \rightarrow T_{u_0}\M$ which is generally \emph{not} the identity map.

This means in particular that parallel transport of a vector along a geodesic of $F_1$ from a point $u_0\in \M$ to $u_1\in \M$ is not equivalent to parallel transport along a geodesic path going from $u_0$ to $u_1$ via $u$; thus even if we parallel transport (in the sense defined by $F_2$) a vector in $W(u)$ along a geodesic of $F_1$ joining $u$ et $u_1$ we cannot be sure that result will lie in $W(u_1)$.

We shall deal with these issues as follows. Assume that $\phi$ satisfies the hypotheses of Lemma~\ref{lemm:step1} and let $u_0\in B(0,R-1)\cap \M$, $u\in B(\M; u_0,\frac{\delta}{2})$ where $ \delta \leq \delta_2$ with $\delta_2 >0$ also given by Lemma~\ref{lemm:step1}.
Choose $w\in W(u)=\textrm{span}~\{T_ue_1,\dots,T_u e_{n+1}\}$, where $(e_1,\dots,e_{n+1})$ is an orthonormal basis of $W\subset T_{u_0}\M$, and assume that $||w|| \leq 1$. Consider the curve: $\sigma(t,u,w)=\exp_u(tw)$, for $t \in [0,t_{max})$ where $t_{max}$ is given by Lemma~\ref{lemm:tmax}. 
Let $\tau \in (0,t_{max})$ and set: \[w_1 = \sigma'(\tau,u, w), \quad u_1=\sigma(\tau,u,w).\]
We begin by an estimate of $||w_1-w||$; by assumption: $\sigma''(t,u,w)=-\frac{|w|^2}{\mu}\sigma(t,u,w)$, thus:
\[ ||\sigma''(t,u,w)|| \leq \frac{|w|^2}{\mu}R.\]
Therefore:
\begin{equation} \label{eq:diffv0vf} ||w_1-w|| \leq \frac{Rt_{\max}}{\mu}.\end{equation}
As we mentioned above $w_1$ may not lie in $W(u_1)$ so we shall consider: \[\displaystyle Pw_1= \sum_{i=1}^{n+1} \langle w_1, T_{u_1}e_i\rangle T_{u_1}e_i,\] the orthogonal projection of $w_1$ onto $W(u_1)$.
In order to estimate $||Pw_1 - w_1||$, we introduce two intermediate vectors. First, let $w_0\in W$ be such that $w=T_uw_0$, we have:
\[ w_0=T_{u}^{-1}w= \sum_{i=1}^{n+1} \langle w,T_ue_i\rangle e_i.\]
Secondly, set $(Pw_1)_0=T_{u_1}^{-1}Pw_1 \in W$, and introduce $v_0=\exp^{-1}_{u_0}(u)$, $v_1=\exp^{-1}_{u_0}(u_1)$. Using Lemma~\ref{lemm:diff_transport} and Equation~\eqref{eq:diffv0vf}:
\[ \begin{aligned} ||Pw_1 -w_1|| &\leq ||Pw_1 -(Pw_1)_0|| + ||(Pw_1)_0-w_0||+||w_0- w|| + ||w - w_1|| \\ &\leq C(R,\mu) |v_1| ||Pw_1|| + ||(Pw_1)_0-w_0|| + C(R,\mu)|v_0|\cdot ||w||  +R\frac{t_{max}}{\mu}\\&\leq C(R,\mu)(|v_1|(1+\frac{Rt_{max}}{\mu})+|v_0|) +R\frac{t_{max}}{\mu} + ||(Pw_1)_0 - w_0||.\end{aligned}\]
Due to holonomy, generally $(Pw_1)_0\neq w_0$, but we can estimate:
\[ \begin{aligned} ||(Pw_1)_0 - w_0||^2 &= \sum_{i=1}^{n+1}(\langle w_1,T_{u_1}e_i)\rangle- \langle w,T_ue_i\rangle)^2\\&=\sum_{i=1}^{n+1}\left( \langle w_1-w,T_{u_1}e_i \rangle +\langle w,T_{u_1}e_i-e_i\rangle + \langle w,e_i-T_ue_i \rangle\right)^2
 \\&\leq (n+1)\left(||w_1-w||+(|v_1|+|v_0|)C(R,\mu)\right)^2 \\&\leq (n+1)\left(\frac{Rt_{max}}{\mu} + (|v_1|+|v_0|)C(R,\mu) \right)^2.\end{aligned} \]
So overall:
\[ \begin{aligned} ||Pw_1 -w_1|| &\leq C(R,\mu)(|v_1|(1+\frac{Rt_{max}}{\mu})+|v_0|) +R\frac{t_{max}}{\mu} +\sqrt{n+1}\left(\frac{Rt_{max}}{\mu}+(|v_1|+|v_0|)C(R,\mu) \right)
\\&\leq C(R,\mu)(3|v_1|+|v_0|)+R\frac{t_{max}}{\mu}+\sqrt{n+1}\left(\frac{Rt_{max}}{\mu} + (|v_1|+|v_0|)C(R,\mu) \right)\\&\leq (3+\sqrt{n+1})C(R,\mu)(|v_0|+|v_1|)+(1+\sqrt{n+1})\frac{Rt_{max}}{\mu}.
\end{aligned}\]
As in the proof of Lemma~\ref{lemm:step1} :
\[ |v_0| \leq \sqrt{\mu} \arccos (1-\frac{\delta}{2\sqrt{\mu}}), |v_1| \leq \sqrt{\mu} \arccos (1-\frac{\delta}{\sqrt{\mu}}).\]
Hence: 
\[ ||Pw_1 - w_1|| \leq 2\sqrt{\mu}(3+\sqrt{n+1})C(R,\mu)\arccos(1-\frac{\delta}{\sqrt{\mu}}) + (1+\sqrt{n+1})\frac{R}{\mu}\frac{\delta}{4}. \]
Now let: 
\begin{eqnarray}\label{def:delta3}
\delta_3 =  \begin{cases} \displaystyle \min\left(\sqrt{\mu}\left(1-\cos\left(\frac{\beta}{\hat{C}}\right) \right), \frac{\beta\mu}{12K(R)R(1+\sqrt{n+1})},\delta_2\right)& \textrm{if $ \displaystyle \frac{\beta}{\hat{C}}<\pi,$}\\ \displaystyle \min\left(\frac{\beta\mu}{12K(R)R(1+\sqrt{n+1})},\delta_2\right) & \textrm{otherwise},\end{cases}
\end{eqnarray}
where $\hat{C}:=96\sqrt{\mu}K(R)(3+\sqrt{n+1})C(R,\mu)$.
With this choice, for any $u\in B(\M; u_0,\delta_3)$, if $\frac{\beta}{\hat{C}}<\pi$ we have setting $K:=K(R)$,
\begin{multline*}
2\sqrt{\mu}(3+\sqrt{n+1})C(R,\mu)\arccos(1-\frac{\delta_3}{\sqrt{\mu}})
\\
\le2\sqrt{\mu}(3+\sqrt{n+1})C(R,\mu)\cdot\frac{\beta}{96\sqrt{\mu}K(R)(3+\sqrt{n+1})C(R,\mu)} = \frac{\beta}{48K},
\end{multline*}
and
$$
(1+\sqrt{n+1})\frac{R}{\mu}\frac{\delta_3}{4}\le (1+\sqrt{n+1})\frac{R}{\mu}\frac{1}{4}\cdot\frac{\beta\mu}{12K(R)R(1+\sqrt{n+1})}=\frac{\beta}{36K}.
$$
Otherwise, we have
\begin{align*}
2\sqrt{\mu}(3+\sqrt{n+1})C(R,\mu)& \arccos(1-\frac{\delta_3}{\sqrt{\mu}})
\\
& \le2\sqrt{\mu}(3+\sqrt{n+1})C(R,\mu)\cdot\frac{\pi}{2}\\
&\le \sqrt{\mu}(3+\sqrt{n+1})C(R,\mu)\cdot\frac{\beta}{96\sqrt{\mu}K(R)(3+\sqrt{n+1})C(R,\mu)}\\
&=\frac{\beta}{96K},
\end{align*}
and $
(1+\sqrt{n+1})\frac{R}{\mu}\frac{\delta_3}{4}\le \frac{\beta}{36K}.$ Thus,
\[||Pw_1 - w_1 || \leq  \frac{\beta}{24K}.\]
It follows then that:
\[ |D^2\phi(u_1)[w_1,w_1] - D^2\phi(u_1)[Pw_1,Pw_1]| \leq ||D^2\phi(u)||\cdot ||Pw_1 -w_1|| \underbrace{||Pw_1+w_1||}_{\leq 6} \leq 6K \frac{\beta}{24K} = \frac{\beta}{4}. \]
Hence:
\[ D^2\phi(u_1)[w_1,w_1] < -\frac{\beta}{2}||Pw_1||^2+\frac{\beta}{4}. \]
We will now assume that $||w||=1$ and avail ourselves of the assumption that $\beta \in (0,1)$,
\[\begin{aligned}  ||Pw_1|| &\geq ||w_1|| - \frac{\beta}{24K} \\& \geq 1- \frac{Rt_{max}}{\mu} - \frac{\beta}{24K} \\&\geq 1- \frac{\beta}{48K}\left( \frac{1}{1+\sqrt{n+1}}+2 \right). \end{aligned}\]
Since $K\geq 1$, $\frac{1}{1+\sqrt{n+1}}+2\leq 3 \leq 3K$, therefore:
\[ ||Pw_1|| \geq 1- \frac{\beta}{16}.\]
Hence, using that $\beta \in (0,1)$:
\[\begin{aligned}-\frac{\beta}{2}||Pw_1||^2+\frac{\beta}{4} &\leq -\frac{\beta}{2}\left(1-\frac{\beta}{16} \right)^2+\frac{\beta}{4} \\ &=-\frac{\beta}{4} +\frac{\beta^2}{16} -\frac{\beta^3}{512} \\&\leq-\frac{\beta}{6}.\end{aligned}\] 
Overall we arrive at:
\[D^2\phi(u_1)[w_1,w_1] < -\frac{\beta}{6}. \]

We can now formulate the following adaptation of \cite[Lemma 3.3]{Fang:1994wz}, key to the proof of Theorem~\ref{M-theorem}.
\begin{lemma}
\label{lemm:3.3}
Let $\phi: E \rightarrow \mathbb{R}$ a $C^2$-functional such that $\phi', \phi''$ are $\alpha$-H\"older continuous ($\alpha \in (0,1]$) on bounded sets. Fix $R>1$ and $u_0 \in B(0,R-1)\cap \M$. Assume that one can find $W\subset T_{u_0}\M$ of finite dimension $n+1$ and $\beta \in (0,1)$ such that for any $w\in W$:
\[ D^2\phi(u_0)[w,w] < -\beta ||w||^2.\]
Let  \[\hat w(u)=\begin{cases}-\frac{P_u\nabla \phi (u)}{||P_u\nabla\phi(u)||}& \textrm{if $ \,\nabla\phi(u)\neq 0$},\\ 0 & \textrm{otherwise},\end{cases}\] where $\nabla \phi$ is the gradient\footnote{The metric dual of the one-form $\textrm{d}\phi|_{\M}$. At each point $u\in \M$, this is simply the orthogonal projection (in the sense of $E$) of the usual gradient $\nabla^E \phi$ in $E$ onto $T_u\M$, i.e. $\nabla \phi(u) = \nabla^E\phi(u) - (\phi'(u)\cdot Gu)\frac{Gu}{||Gu||^2} \in T_u\M={Gu}^\perp.$ } of $\phi$ and $P_u$ denotes the orthogonal projection (in $T_u\M$) onto $T_u W$. \smallskip

\noindent
 Then, if $u\in B(\M; u_0,\frac{\delta_3}{2})$, where $\delta_3 >0$ satisfies \eqref{def:delta3}, and $t_{max}=\min(\frac{2\mu}{R}, \frac{\delta_3}{4})$ either:
\begin{enumerate}
 \item $\hat w(u)\neq 0$ and for any $t\in (0,t_{max})$, we have:
 \[ \phi(\exp_u(t \hat w(u)))< \phi(u) -\frac{\beta}{12} t^2;\] 
 or 
 \item $\hat w(u)=0$, and for any $w\in T_u W$, $||w||=1$ and any $t\in (0,t_{max})$,
 \[\phi(\exp_u(tw))<\phi(u) -\frac{\beta}{12}t^2.\]
 Furthermore, in this case, for any $t_0\in (0,t_{max})$, one can find $0<\delta_u < \frac{\delta_3}{2}-|u-u_0|$ such that for any  $z\in \overline{B}(u,\delta_u)\cap \M$, $w \in T_z W$ and $t \in [t_0,t_{max})$ we have:
 \[\phi(\exp_z(tw))<\phi(z) -\frac{\beta}{24}t^2.\]
 \end{enumerate}
\end{lemma}
\begin{proof}
For every $w\in T_u W$, and any $t\in (0,t_{max})$ the Taylor-Lagrange theorem guarantees that there is $\tau \in (0,t)$:
\[\phi(\exp_u(tw)) =\phi(u) + t \langle \nabla \phi(u), w\rangle + \frac{t^2}{2}D^2\phi(\exp_u(\tau w))\left[\left.\frac{d}{dt}(\exp_u(tw))\right|_{t=\tau},\left.\frac{d}{dt}(\exp_u(tw))\right|_{t=\tau}\right].\]

\begin{enumerate}
\item If $\hat w(u)\neq 0$, substitute in the above $w=\hat w(u)$ and notice that:
\[ \langle \nabla \phi(u), \hat w(u) \rangle= -1.\]
Since $t>0$ it follows that:
\[ \phi(\exp_u(t\hat w(u))) <\phi(u) + \frac{t^2}{2}D^2\phi(\exp_u(\tau \hat w(u)))\left[\left.\frac{d}{dt}(\exp_u(t\hat w(u)))\right|_{t=\tau},\left.\frac{d}{dt}(\exp_u(t\hat w(u)))\right|_{t=\tau}\right]. \]
\item If $\hat w(u) =0$, then take any $w \in T_u W$ with $||w||=1$, and notice that $\langle \nabla \phi(u), w \rangle = 0$, therefore:
\[ \phi(\exp_u(tw)) =\phi(u) + \frac{t^2}{2}D^2\phi(\exp_u(\tau w))\left[\left.\frac{d}{dt}(\exp_u(tw))\right|_{t=\tau},\left.\frac{d}{dt}(\exp_u(tw))\right|_{t=\tau}\right]. \]
\end{enumerate}
In order to conclude the proof in both cases we apply the discussion preceding the lemma which implies that, for any $u\in B(\M; u_0,\frac{\delta_3}{2})$, $w\in T_uW$, $||w||=1$ and any $\tau \in (0,t_{max})$: 
\begin{equation} \label{eq:EstimD2} D^2\phi(\exp_u(\tau w))\left[\left.\frac{d}{dt}(\exp_u(tw))\right|_{t=\tau},\left.\frac{d}{dt}(\exp_u(tw))\right|_{t=\tau}\right]< -\frac{\beta}{6}.\end{equation}
The first two points ensue.

For the final point, by continuity of $\phi'$ at $u$, choose $\delta_u < \frac{\delta_3}{2} -|u-u_0|$ such that
for any  $z\in \overline{B}(u,\delta_u)\cap \M$:
\begin{equation}\label{eq:E_1} ||\phi'(u)- \phi'(z)||<\frac{\beta}{48}t_0.\end{equation}
We must however take into account the fact that $T_uW \neq T_z W$ so we need to approximate $w_z \in T_z W$ by $w_u \in T_u W$. To this end, observe that for every $u$, $T_u$ can be viewed as a continuous map $T_u : T_{u_0}\M \rightarrow T_u \M \subset E$.
The restriction of $T_u$ to $W$, denoted by $T_u|_W$, can be thought of as an element of $L(W,E)$. We claim that the map $u\mapsto T_u|_W$ is continuous.

For this we shall use the fact that $W$ is finite dimensional; choose an orthonormal basis $(e_1,\dots, e_{n+1})$ of $W$, if for each $i\in \{1,\dots, n+1\}$: $u\mapsto T_ue_i$ is continuous,  it will follow that $u\mapsto T_u|_W$ is continuous. Indeed, let $\varepsilon>0$ and observe that $||\cdot ||$ and\footnote{$(e_i^*)$ is the dual basis to $(e_i)$} $|w|_\infty=\max_{i\in \{1,\dots,n+1\}} |e^*_i(w)|$ define equivalent norms on $W$. Therefore, if
$u \in B(\M; u_0, \frac{\delta_3}{2})$ is fixed and $\tilde{\delta}$ is chosen small enough such that for every $i\in\{1,\dots, n+1\}$ $z \in \overline{B}(u,\tilde{\delta})\cap \M$: \[||T_u e_i -T_z e_i || \leq \frac{\varepsilon}{n+1}, \]
we have
\begin{equation*}
||T_u w - T_z w ||\leq\varepsilon|w|_\infty \leq C \varepsilon ||w||,~~~ \textrm{for all $w\in W$,}
\end{equation*}
where $C>0$ is a constant independent of $w\in W$.
Therefore $u\mapsto T_u|_W$ is continuous. It remains to justify that the map $u\mapsto T_ue_i$ is continuous for each $i\in\{1,\dots,n+1\}$. 
Fix such an $i$, and recall that $T_ue_i$ is defined to be $\varphi(1)$, where $\varphi$ is the solution to the Cauchy problem:
\[\varphi'(t) + (\varphi(t),\frac{d}{dt}(\exp_{u_0}(tv))\frac{G\exp_{u_0}(tv)}{||G\exp_{u_0}(tv)||^2}=0, \quad \varphi(0)=e_i,\]
with $v=\exp^{-1}_{u_0}u$. Denoting the differential of the map $v\mapsto \exp_{u_0}v$ at a point $v$ by ${\exp_{u_0}}_{*v}$, the above can be rewritten as a system:
\[ \begin{cases} \varphi'(t)+(\varphi(t),{\exp_{u_0}}_{*t\psi(t)}\cdot \psi(t)) \frac{G\exp_{u_0}(t\psi(t))}{||G\exp_{u_0}(t\psi(t))||^2}=0,\\ \psi'(t)=0,\end{cases}\]
with initial conditions: $\psi(0)=v$, $\varphi(0)=e_i$. 

By continuity with respect to initial conditions, the solution, and a fortiori,  $\varphi(1)$, depend continuously on $v$, however $v=\exp_{u_0}^{-1}(u)$, so $\varphi(1)$ and thus $T_ue_i$ depend continuously on $u$.
This proves our claim. Consequently, reducing further, if necessary, $\delta_u$ we can assume that for any $z\in \overline{B}(u,\delta_u)\cap \M$:
 \begin{equation}\label{eq:E2} ||T_u|_W-T_z|_W||<\frac{\beta \, t_0}{48K(R)},\end{equation}
where $K(R)$ is defined in \eqref{def:K}.
Recall now that, for any $w\in T_z W$, $||w||=1$, $z\in \overline{B}(u,\delta_u)\cap \M \subset B(\M;u_0,\frac{\delta_3}{2})$, $t\in [t_0,t_{max})$ one can find $\tau \in (0,t)$ such that:
 \[\phi(\exp_z(tw)) = \phi(z) +t\phi'(z)\cdot w +\frac{t^2}{2}D^2\phi(\exp_z(\tau w))\left[\left.\frac{d}{dt}(\exp_z(tw))\right|_{t=\tau},\left.\frac{d}{dt}(\exp_z(tw))\right|_{t=\tau}\right].\]

If we write $w_0=T_z^{-1}w$, then $\phi'(u)\cdot T_u w_0=0$, and, using Lemma \ref{para tran}, Equations \eqref{def:K}, \eqref{eq:E_1} and \eqref{eq:E2}:
\[\begin{aligned}|\phi'(z)\cdot w|&=|\phi'(z)\cdot T_z w_0-\phi'(u)\cdot T_u w_0|\\&\leq |\phi'(z)\cdot T_z w_0-\phi'(u)\cdot T_z w_0| + |\phi'(u)\cdot T_z w_0-\phi'(u)\cdot T_uw_0| \\&\leq  ||\phi'(z)-\phi'(u)||\cdot ||T_z w_0||+||\phi'(u)||\cdot||T_z w_0 - T_u w_0|| \\
&\le\frac{\beta}{48}t_0+K(R)\cdot\frac{\beta t_0}{48K(R)}\\
&= \frac{\beta}{24}t_0. \end{aligned}\]
Since $z\in \overline{B}(u_0,\frac{\delta_3}{2})\cap \M$, we can estimate the second order term using Equation~\eqref{eq:EstimD2} and thus,  for $t\in[t_0, t_{\max})$:
\[\begin{aligned} \phi(\exp_z(tw)) &< \phi(z) +t\frac{\beta}{24}t_0 -\frac{\beta}{12}t^2 \\ &<\phi(z) - \frac{\beta}{24} t^2.\end{aligned}\]
\end{proof}

%
%
%
%
%
\begin{lemma}\label{lemm:3.4-FG}
Assume that the hypotheses of Lemma~\ref{lemm:3.3} are satisfied and let $f$ be a continuous map from a closed subset $D\subset \mathbb{R}^n$ into $\M$. Suppose that $K_1 \subset D$ is a compact subset of $D$ such that $f(K_1)\subset B(\M;u_0,\frac{\delta_3}{2})$. Then for sufficiently small $\nu>0$ and $t_0\in (0,t_{max})$  there is a continuous map $\eta: [0,t_{max}] \times D \rightarrow \M$ that satisfies:
\begin{enumerate}
\item $\eta(t,x)= f(x)$, if $(t,x) \in (\{0\}\times D) \cup ([0,t_{max}]\times D\setminus N_\nu(K_1))$;
\item $\phi(\eta(t,x)) \leq \phi(f(x))$ if $(t,x) \in [0,t_{max})\times D$;
\item $\phi(\eta(t,x)) < \phi(f(x)) -\frac{\beta}{24}t^2$ if $(t,x)\in[t_0,t_{max})\times K_1$;
\item $||\eta(t,x)-f(x)|| \leq 3t $ for all $(t,x)\in [0,t_{max}]\times D$.
\end{enumerate}
\end{lemma}
\begin{proof}
Let $N_\nu(K_1)$ denote the $\nu$-neighbourhood of $K_1$ in $\mathbb{R}^n$, namely,
\[ N_\nu(K_1)=\{x \in \mathbb{R}^n \,|\, d(x,K_1)<\nu \}.\]
Assume that $\nu >0$ is small enough so that, $f(\overline{N_\nu(K_1)}\cap D)\subset B(\M; u_0,\frac{\delta_3}{2})$ and define:
\[T=\{x \in \overline{N_\nu(K_1)}\cap D \,|\, \hat{w}(f(x))=0\}.\]
By the last point in Lemma~\ref{lemm:3.3}, for every $y\in T$ one can find $\nu^y >0$ such that for any $z\in B(y,\nu^y)\cap D$, and any $w\in T_{f(z)}W$, $||w||=1$ we have the inequality:
\[\phi(\exp_{f(z)}(tw)) < \phi(f(z)) - \frac{\beta}{24}t^2,\quad  t\in [t_0, t_{max}).\]
Put $O=\cup_{y\in T} B(y,\frac{\nu^y}{2})$ and let $g:\mathbb{R}^n \rightarrow [0,1]$ be a continuous function such that: 
\[ g(x)=\begin{cases} 1 & x\in K_1 ,\\ 0 & x\in \mathbb{R}^n\setminus N_{\nu}(K_1). \end{cases}\]

Next, choose an orthonormal basis $(e_1,\dots,e_n)$ of $W \subset T_{u_0}\M$ and let $f_1: \overline{N_\nu(K_1)}\cap D\setminus O \rightarrow S^{n} \subset \mathbb{R}^{n+1}$ be defined by:
\[ f_1(x)=(\langle T_{f(x)} e_1, w(f(x))\rangle, \dots, \langle T_{f(x)} e_{n+1}, w(f(x))\rangle).   \]

As in the proof of Lemma~\ref{lemm:3.3}, the continuity of the map: $u\mapsto T_u|_W$ guarantees that $f_1$ is itself continuous.  
Since $\overline{N_\nu(K_1)}\cap D\setminus O$ is closed in $\mathbb{R}^n$ one can extend $f_1$ to a continuous map $f_2: \mathbb{R}^n \rightarrow S^n$ by Lemma~\ref{Corollary 2.4-FG}.

Define now a continuous map $f_3$ on $D$:
\[f_3(x) = \sum_{i=1}^{n+1} E_i^*(f_2(x))T_{f(x)}e_i, \in T_{f(x)}W, \]
where we denote $(E_i^*)$ the dual basis of the canonical basis $(E_1,\dots, E_{n+1})$ of $\mathbb{R}^{n+1}$,
then set \[\eta(t,x)=\exp_{f(x)}(t g(x)f_3(x)), x \in D, t \in [0,t_{max}].\]
It is straight forward to check that $\eta$ satisfies the required conditions. The last point follows from Equation~\eqref{eq:diffv0vf} and the mean value theorem.
\end{proof}

\begin{lemma}\label{Lemma3.5-FG}
Let $\phi: E \rightarrow \mathbb{R}$ a $C^2$-functional such that $\phi', \phi''$ are $\alpha$-Hölder continuous ($\alpha \in (0,1]$) on bounded sets and let $f$ be a continuous map from a compact subset $D$ of $ \, \mathbb{R}^n$ into $\M$.
Suppose $K_2$ is a compact subset of $D$ with the following properties:  
\begin{itemize}
\item There exists $R >1$ such that $f(K_2) \subset B(0, R-1)\cap \M$.
\item There exists a constant $\beta >0$ such, that for all $y \in K_2$, there is a subspace $W_y$, of $T_{f(y)}\M$ with $\dim \, W_y \geq n+1$ so that 
\begin{equation}\label{3.2-FG-l}
D^2 \phi(f(y))[w,w]< - \beta ||w||^2, \quad \mbox{for all } \, w \in W_y.
\end{equation}
\end{itemize}

Then for any $0 < \delta \leq \delta_3$ where $\delta_3 >0$ satisfies \eqref{def:delta3} and $\nu >0$ there is a continuous map  $\hat{f} : D \rightarrow \M$ such that if $N:= N(n)$ is the number given in 
Lemma \ref{Lemma 2.5-FG}, we have : \smallskip

\begin{itemize}
\item[(i)] $\hat{f}(x)= f(x)$ for $u \in D \backslash N_{\nu}(K_2)$; \smallskip
\item[(ii)] $\phi(\hat{f}(x)) \leq \phi (f(x))$  for all $x \in D$; \smallskip
\item[(iii)] If $x \in K_2$, then $\phi(\hat{f}(x)) <  \displaystyle \phi (f(x)) - \frac{\beta \, \delta^2}{864 N^2};$ \smallskip
\item[(iv)] $ ||\hat{f}(x) - f(x) || \leq \frac{\delta}{2}$ for all $x \in D$.
\end{itemize}
\end{lemma}

\begin{proof}
Let $0 < \delta \leq \delta_3$ be fixed. Since $f$ is uniformly continuous on  $D \subset \R^n$, there exists $0 < \varepsilon_N < \frac{\nu}{4}$ such that
\begin{equation}\label{0.8}
||f(x) - f(y)|| \leq \frac{\delta}{8N}, \quad \mbox{for all  } (x,y) \in D \times D \, \mbox{ with } \, d(x,y) \leq 2 \varepsilon_N.
\end{equation}

Using Lemma \ref{Lemma 2.5-FG} with $C = K_2$ and $\varepsilon = 4 \varepsilon_N$ we deduce that there exists a finite number of distinct points $\{x_i, 1 \leq i \leq m \}$ of $\R^n$ such that $\cup_{i=1}^{m}B(x_i, \varepsilon_N)$ covers $K_2$ and such that any intersection of $N$ distinct $\overline{B(x_i, 2\varepsilon_N)}$ is empty. 

We may certainly assume that $B(x_i,\varepsilon_N)\cap K\neq\emptyset$ for each $i \in \{1, \cdots, m\}$ and choose $y_i \in B(x_i, \varepsilon_N) \cap K_2$. For convenience, we set $B_{y_i} = B(x_i, \varepsilon_N)$. Observe, from \eqref{0.8} that
$$f(x) \in B \Big(\M; f(y_i), \frac{\delta}{8N} \Big), \quad x \in \overline{B_{y_i}} \cap D.$$ 
Choose  $0 < \tau < \va_N$ small enough such that $\cup_{i=1}^m N_{\tau}(B_{y_i}) \subset N_{\nu}(K_2)$ and for $i \in 1, \cdots, m$
\begin{equation*}
f(x) \in B \Big(\M; f(y_i), \frac{\delta}{4N} \Big), \quad x \in \overline{N}_{\tau}(B_{y_i}) \cap D.
\end{equation*}
Note also that, since $\tau + \varepsilon_N < 2 \varepsilon_N$, any intersection of $N$  distinct sets $\overline{N_{\tau}}(B_{y_i})$ is empty. \\
We shall now define by induction, continuous functions $f_0, f_1, \cdots, f_m : D \rightarrow \M$ such that for all $1 \leq i \leq m$ we have that
\begin{equation}\label{3.3-l}
\phi(f_i(x)) < \phi(f_{i-1}(x)) - \frac{\beta \, \delta^2}{864 N^2} \quad \mbox{if } x \in \overline{B}_{y_i} \cap D,
\end{equation}
\begin{equation}\label{3.4-l}
\phi (f_i(x)) \leq \phi (f_{i-1}(x)) \quad \mbox{if } x \in  D,
\end{equation}
and
\begin{equation}\label{3.5-l}
||f_i(x) - f_{i-1}(x)|| \leq
\begin{cases}
0  \quad \mbox{if } x \in D \backslash N_{\tau}(\overline{B}_{y_i}),\\
\frac{\delta}{2 N}\quad \mbox{if } x \in N_{\tau}(\overline{B}_{y_i}) \cap D.
\end{cases}
\end{equation}
Let $f_0= f$ and suppose that $f_0, f_1, \cdot, f_k$ are well-defined and satisfy inequalities \eqref{3.3-l}, \eqref{3.4-l} and \eqref{3.5-l} for $k <m.$ Clearly
$$||f_i(x) - f(x) || \leq \frac{i \delta}{2 N} \quad \mbox{if } x \in \bigcap_{j=1}^i N_{\tau}(\overline{B}_{y_i}) \cap D.$$
Since any intersection of $N$ distinct sets $N_{\tau}(B_{y_i})$ is empty, we have that
$$ ||f_k(x) - f(x) || \leq \frac{\delta(N-1)}{2  N} \quad \mbox{if } x \in  D.$$
Since $f : \overline{B}_{y_{k+1}} \cap D \rightarrow B(\M; f(y_{k+1}), \frac{\delta}{8  N})$, we see that $f_k$ 
maps $\overline{B}_{y_{k+1}} \cap D$ into 
$B(\M; f(y_{k+1}), \frac{\delta (N-1)}{2 N}+ \frac{\delta}{8 N}) \subset B(\M; f(y_{k+1}), \frac{\delta}{2 }).$

By assumption \eqref{3.2-FG-l}, there is some subspace $W_{y_{k+1}}$ of $E$ with $\dim W_{y_{k+1}} \geq n+1$ such that for any $w \in W_{y_{k+1}},$ with $||w|| =1,$ we have that 
$D^2 \phi(f(y_{k+1}))[w,w]< - \beta$.
Hence, we may apply Lemma \ref{lemm:3.4-FG} with $f_k$ and any $t_0 \in (0, t_{max})$ to obtain a continuous deformation $\eta(t,x)$ satisfying the conclusion of that lemma. Define now $f_{k+1}(x) = \eta (\frac{\delta}{6 N},x)$ to get a continuous function  $f_{k+1}: D \rightarrow \M$ satisfying
$$ \phi(f_{k+1}(x)) < \phi(f_{k}(x)) - \frac{\beta \, \delta^2}{864 N^2} \quad \mbox{for } x \in \overline{B}_{y_{k+1}} \cap D,$$
\begin{equation*}
\phi (f_{k+1}(x)) \leq \phi (f(x)) \quad \mbox{for } x \in  D,
\end{equation*}
and
\begin{equation*}
||f_{k+1}(x) - f_{k}(x)|| \leq
\begin{cases}
0  \quad \mbox{if } x \in D \backslash N_{\tau}(\overline{B}_{y_{k+1}}),\\
\frac{\delta}{2 N}\quad \mbox{if } x \in N_{\tau}(\overline{B}_{y_{k+1}}) \cap D.
\end{cases}
\end{equation*}
By induction we see that $f_0, \cdots, f_m$ are well-defined. Clearly $\hat{f} = f_m$ verifies the claims of the lemma.
\end{proof}

Finally, we shall need the following lemma which follows directly from \cite[Lemma 5.15]{Willem} used with $S = f(K_3)$ and $V = \M$.
\begin{lemma}\label{Lemma3.6-FG}
Let $\phi$ be a $C^1$ functional on $E$ and let $f$ be a continuous map from a closed subset $D$ of $\, \mathbb{R}^n$ into $\M$. Let $\tilde{c}$, $\tilde{\va}$, $\tilde{\mu} >0$
 be three constants.
 Suppose $K_3$ is a compact subset of $D$ satisfying
$$ \tilde{c} - \tilde{\va} \leq \phi(f(x)) \leq \tilde{c} + \tilde{\va}, \quad \mbox{for } x \in K_3.$$
Assume that, for all $x \in K_3$,
$$||\phi'|_{\M}(u) || \geq \frac{8 \tilde{\va}}{\tilde{\mu}}, \quad \mbox{for } u \in B(\M; f(x), 2 \tilde{\mu}), $$
then there is a continuous map  $\hat{f} : D \rightarrow \M$ such that 
\smallskip

\begin{itemize}
\item[(i)] $\hat{f}(x)= f(x)$ if $ \,\phi(f(x)) \leq \tilde{c} - 2 \tilde{\va}$; \smallskip
\item[(ii)] $\phi(\hat{f}(x))) \leq \phi (f(x))$  for all $x \in D$; \smallskip
\item[(iii)] If $x \in K_3$, then $\phi(\hat{f}(x))\leq  \tilde{c} - \tilde{\va}.$
\end{itemize}
\end{lemma}


We now have all the ingredients to prove the main result of this section.

\begin{proof}[Proof of Theorem \ref{M-theorem}]
 
Suppose $\max_{u \in A}\phi(u) \leq c + \va$ where $A$ is a set in $\mathcal{F}$ satisfying \eqref{top-bounded}. There exists a continuous function $f$ from $D \subset \mathbb{R}^n$ into $\M$, which is equal to $\sigma$ on $D_0$ and such that $A = f(D)$. 

Let $\delta_3 >0$ satisfy the equality in \eqref{def:delta3}. Observing that it is not restrictive to assume that $\alpha \leq \frac{1}{2}$ in condition \eqref{Holder} we can find a constant $\gamma >0$ such that $\gamma \beta^{\frac{1}{\alpha}} \leq \delta_3$ for any $\beta >0$ small enough. Now let
\begin{equation}\label{value-delta}
\delta = \delta(\va) = \frac{1}{2}\min\Big( \frac{\va^{\alpha_1}}{M}, \gamma \, \va^{\frac{\alpha_1}{\alpha}}\Big)
\end{equation}
and observe that $\delta \leq \delta_3$ when $\beta >0$ is given by $\beta = \va^{\alpha_1}$.

Consider the closed set
\begin{equation}\label{def-K}
K = \{ x \in D \,|\, \phi(f(x)) \geq c- \varepsilon \}.
\end{equation}
Since $c - \max_B \phi >0$, taking $\varepsilon >0$ sufficiently small, we can assume that $K$ is a compact subset of $D \backslash D_0$  and that $c - 2 \va > \max_B \phi$. \\
Now suppose that the conclusion of the theorem does not hold, then for all $x \in K$, we have either $||\phi'|_{\M}(f(x))|| > 3\va^{\alpha_1}$ or there is a subspace $W_x$ of $T_{f(x)}\M$ with $\dim W_x \geq n+1$ such that for all $w \in W_x$, we have $D^2 \phi(f(x))[w,w] < -  \varepsilon^{\alpha_1} ||w||^2$. \medskip 

In view of the assumption \ref{hol2} on $\phi'$ and of the definition of $\delta >0$ given in \eqref{value-delta}, setting $\hat{\delta}= 2 \delta$ we deduce that for all $x \in K$ such that  $||\phi'|_{\M}(f(x))|| > 3\va^{\alpha_1}$ we have $||\phi'|_{\M}(u)|| > \va^{\alpha_1}$ for all $u \in B(f(x), \hat{\delta})$.
Let 
$$T_1= \{x \in K \,|\, ||\phi'|_{\M}(u)|| > \varepsilon^{\alpha_1}, \mbox{ for all } u \in B(\M;f(x),\hat{\delta})\}$$ and $T_2 = K \backslash T_1$.
  Note that $\overline{T}_1$, $\overline{T}_2$ are compact, and $K = \overline{T}_1 \cup \overline{T}_2$. Now apply Lemma \ref{Lemma3.5-FG} with $K_2 = \overline{T}_2$, $\beta = \va^{\alpha_1}$  and $\nu = \frac{1}{2}\dist(D_0, K) >0$ to obtain a continuous map  $g : D \rightarrow \M$ such that
\begin{equation}\label{3.13}
g(x) = f(x) \quad \mbox{for } x \in D \backslash N_{\nu}(\overline{T}_2) \quad \mbox{and} \quad \phi(g(x)) \leq \phi (f(x)), \quad \mbox{for } x \in D,
\end{equation}
\begin{equation}\label{3.14}
\phi(g(x)) \leq \phi (f(x)) - \frac{\varepsilon^{\alpha_1} \delta^2}{864  N^2} , \quad \mbox{for } x \in \overline{T}_2,
\end{equation}
\begin{equation}\label{3.15}
||g(x) - f(x) || \leq \frac{\delta}{2}, \quad \mbox{for } x \in N_{\nu}(\overline{T}_2) \cap D.
\end{equation}
Observe that Inequality \eqref{3.15} yields for $x \in \overline{T}_1,$  $B(\M; g(x), \delta)  \subset B(\M;f(x), \hat{\delta})$. Now we apply Lemma \ref{Lemma3.6-FG} with $\tilde{c} =c$, $\tilde{\va} = \va$ and $\tilde{\mu} = \frac{\delta}{2}$. Observing that since $\alpha_1 < \frac{\alpha}{\alpha +1}$ we have that
$$\va^{\alpha_1} \geq \frac{16 \va}{\delta}, \quad \mbox{if } \va >0, \mbox{ is sufficiently small.}$$
we deduce that there exists  a continuous map  $\hat{f} : D \rightarrow \M$ such that
\begin{equation}\label{3.17}
\hat{f}(x) = g(x) \quad \mbox{if } \phi(g(x)) \leq c- 2 \va  \quad \mbox{and} \quad \phi(\hat{f}(x)) \leq \phi(f(x)), \quad \mbox{for} \quad x \in D,
\end{equation}
\begin{equation}\label{3.18}
\phi(\hat{f}(x)) \leq c - \va, \quad \mbox{for} \quad x \in \overline{T}_1,
\end{equation}
Note that $\hat{f}(D) \in \mathcal{F}$. In view of \eqref{3.14}, \eqref{3.17}, \eqref{3.18} and of the definition of $K$ given in \eqref{def-K}  to get a contradiction we just need to show that
\begin{equation}\label{key-estimate}
\frac{\va^{\alpha_1}\delta^2}{864  N^2} > 2 \varepsilon, \quad \mbox{for} \quad x \in \overline{T}_2.
\end{equation}
But, for $\va >0$ small enough,
\begin{equation*}
 \delta^2 =  \frac{1}{4}\gamma^2 \va^{\frac{2 \alpha_1}{\alpha}}
\end{equation*} 
and since $\alpha_1 < \frac{\alpha}{\alpha +2}$, Equation \eqref{key-estimate} will hold if 
$\va >0$ is small enough. \smallskip

Summarising, we managed to construct, assuming that the conclusion of the theorem does not hold, a path in $\mathcal{F}$ for which the maximum of $\phi$ on this path is strictly below the value of $c$. This contradiction ends the proof.
\end{proof}

\section{Proof of Theorem \ref{thm: application}}\label{Theorem-application}

This section is devoted to the

\begin{proof}[Proof of Theorem \ref{thm: application}]
Since, for each $u \in E$, the function $\rho \mapsto \Phi_\rho(u)$ is non-increasing, the function $\rho \mapsto c_\rho$ is non-increasing as well. Therefore, its derivative $c_\rho'$ is well defined for almost every $\rho \in I$. We show that the existence of $c_\rho'$ ensures that of the desired Palais-Smale sequence. Let then $\rho \in I^\circ$ (the interior of $I$) be such that $c_\rho'$ exists, and let $\{\rho_n\} \subset I$ be a monotone increasing sequence converging to $\rho$.

\medskip

\emph{Step 1)} There exist $\{A_n\} \subset \mathcal{F}$ and $K=K(c_\rho')>0$ such that, writing $A_n = f_n(D)$ we have:
\begin{itemize}
\item[(i)] $\|f_n(x)\| \le K$ whenever
\begin{equation}\label{2.1}
\Phi_\rho(f_n(x)) \ge c_\rho - (2-c_\rho')(\rho-\rho_n),~ x\in D;
\end{equation}
\item[(ii)] $\displaystyle
 \max_{x \in D} \Phi_\rho(f_n(x)) \le c_\rho + (2-c_\rho')(\rho-\rho_n).$
\end{itemize}
Let $\{A_n\} \subset \mathcal{F}$ be an arbitrary sequence such that,  writing $A_n = f_n(D)$ we have
\begin{equation}\label{2.2}
\max_{x \in D}\Phi_{\rho_n}(f_n(x)) \leq c_{\rho_n}+ (\rho - \rho_n).
\end{equation}
We shall prove that, for $n \in \N$ sufficiently large, $\{A_n\} \subset \mathcal{F}$ satisfies the desired conditions.
When $f_n(x)$ satisfies \eqref{2.1}, we
have
$$\frac{\Phi_{\rho_n}(f_n(x))- \Phi_{\rho}(f_n(x))}{\rho - \rho_n } \leq \frac{c_{\rho_n} + (\rho - \rho_n) - c_{\rho} + (2 - c'_{\rho})(\rho- \rho_n)}{\rho- \rho_n} = \frac{c_{\rho_n} - c_{\rho}}{\rho - \rho_n} + (3 - c'_{\rho}).$$
Since $c'_{\rho}$ exists, there is a $n(\rho) \in \N$ such that, for all $n \geq n(\rho)$,
\begin{equation}\label{2.0}
\frac{c_{\rho_n}- c_{\rho}}{\rho - \rho_n} \leq - c'_{\rho} +1.
\end{equation}
Consequently, for all $n \geq n(\rho)$,
\begin{equation}\label{1}
B(f_n(x)) = \frac{\Phi_{\rho_n}(f_n(x))- \Phi_{\rho}(f_n(x))}{\rho - \rho_n } \leq - 2c'_{\rho} +4.
\end{equation}
Moreover, 
\begin{equation}\label{2}
A(f_n(x)) = \Phi(f_n(x)) + \rho_n B(f_n(x)) \leq c_{\rho_n} + (\rho - \rho_n) + \rho_n (-2 c'_{\rho} +4).
\end{equation}
(i) now follows from Assumption \eqref{hp coer} and the fact that $\{A(f_n(x))\}$ and $\{B(f_n(x))\}$ are bounded. To prove (ii), observe from \eqref{2.0} that for all $n \geq n(\rho)$,
\begin{equation}\label{2.4}
c_{\rho_n} \leq c_{\rho} + (- c'_{\rho}+1) (\rho - \rho_n)
\end{equation}
and thus, using \eqref{2.2} and \eqref{2.4}, we get
$$ \Phi_{\rho}(f_n(x)) \leq \Phi_{\rho_n}(f_n(x)) \leq c_{\rho_n}+ (\rho- \rho_n) \leq c_{\rho} + (- c'_{\rho}+2) (\rho - \rho_n).$$
Thus (ii) also holds.
\medskip

\emph{Step 2)} Let $K=K(c_\rho')$ be the constant found in Step 1, and let $M \ge 1$ be the constant of the assumption \ref{hol2}
 corresponding to the value $R=K+1$. 
Finally, let $\eps_n:= (2-c_\rho')(\rho-\rho_n) \to 0^+$, (where $\{\rho_n\}$ is the sequence introduced at Step 1). For each  $n \in \N$ large enough, Theorem \ref{M-theorem} guarantees the existence of 
 $u_n \in \M$ with the following properties:
\begin{itemize}
\item[(a)] $c_\rho-\eps_n \le \Phi_\rho(u_n) \le c_\rho+\eps_n$;
\item[(b)] $ \|\Phi'_\rho|_{\M}(u_n)\|_* \le  3\eps_n^{\alpha_1}$;
\item[(c)] $u_n \in A_n$;
\item[(d)] $\tilde m_{\eps_n^{\alpha_1}}(u_n) \le d$.
\end{itemize}

\emph{Step 3) Conclusion of the proof.}
Recalling that $\eps_n =(2-c_\rho')(\rho-\rho_n) \to 0^+$, we claim that $\{u_n\} \subset \M$ with $u_n \in \M$ satisfying (a)-(d) of Step 2 is the Palais-Smale sequence we are looking for in Theorem \ref{thm: application}.\\
From (a) and (b) of Step 2, it is a Palais-Smale sequence at the level $c_\rho$. From (a) and (c) we also have  
\[
u_n \in A_n \cap \{\Phi_\rho \ge c_\rho-(2-c_\rho')(\rho-\rho_n)\},
\]
so that $\|u_n \| \le K$, thanks to Step 1; that is, $\{u_n\} \subset \M$ is bounded.  Finally, we observe that, since $\va_n^{\alpha_1} \to 0^+$, (d) of Step 2 implies (iv) of Theorem \ref{thm: application}. At this point the proof of the theorem is completed. 
\end{proof}

\begin{proof}[Proof of Remark \ref{space-of-functions}]
The remark follows directly from the observation that it is possible to replace in Step 1 of the proof of Theorem \ref{thm: application}, $A_n \in \mathcal{F}$ by $|A_n|_* \in \Gamma$. Since $u_n \in |A_n|_* $, this gives the additional property $u_n \ge 0$.
\end{proof}

Theorem \ref{thm: application} admits variants in which some of the assumptions can be relaxed. The following one is motivated by a remark of D. Ruiz. 
\begin{theorem}\label{thm: application-M}
Let $I \subset (0,+\infty)$ be an interval and consider a family of $C^1$ functionals $\Phi_\rho: E \to \mathbb{R}$ of the form
\[
\Phi_\rho(u) = A(u) -\rho B(u), \qquad \rho \in I,
\]
where $B(u) \ge 0$ for every $u \in E$.
Let $\mathcal{F}$ be a homotopic family of $\M$ of dimension at most $d$ with boundary $B$ (independent of $\rho$) such that 
\begin{equation}\label{R-mp geom}
c_{\rho}:= \inf_{A \in \mathcal{F}} \max_{u \in A}\Phi_{\rho}(u) > \max_B \Phi_{\rho},\quad \forall \rho \in I.
\end{equation}
Then, for almost every $\rho \in I$, there exists a sequence $\{u_n\} \subset \M$  such that, as $n \to + \infty$,
\begin{itemize}
\item[(i)] $\Phi_\rho(u_n) \to c_\rho$;
\item[(ii)] $||\Phi'_\rho|_{\M}(u_n)||_* \to 0$; 
\item[(iii)] $\{A(u_n)\}$ and $\{B(u_n)\}$ are bounded.
\end{itemize}
Assuming in addition that the set $G:= \{ u \in \M \,|\,  A(u) \leq K_1, B(u) \leq K_2 \}$ is bounded for any $K_1 >0$, $K_2 >0$, then also $\{u_n\}$ is bounded in $E$. \\
If now $\Phi_\rho: E \to \mathbb{R}$  is of class $C^2$ and $\Phi_\rho'$, $\Phi_\rho''$ are $\alpha$-Hölder continuous on bounded sets in the sense of Definition~\ref{Holder continuous}~(\ref{hol2}) for some $\alpha \in (0,1]$, then there exists a sequence $\zeta_n \to 0^+$ such that, as $n \to + \infty$,
\begin{itemize}
\item[(iv)] $\tilde m_{\zeta_n}(u_n) \le d$.
\end{itemize}
\end{theorem}
\begin{proof}
First by a direct adaptation of Step 1) of the proof of Theorem \ref{thm: application} we obtain, see \eqref{1} and \eqref{2}, that there exist $\{A_n\} \subset \mathcal{F}$, $A_\rho >0$ and $B_{\rho}>0$ such that, writing $A_n = f_n(D)$ we have
\begin{itemize}
\item[(i)] $A(f_n(x)) \leq A_{\rho}$ and $B(f_n(x)) \leq B_{\rho}$ whenever
\begin{equation*}
\Phi_\rho(f_n(x)) \ge c_\rho - (2-c_\rho')(\rho-\rho_n);
\end{equation*}
\item[(ii)] $\displaystyle
 \max_{x \in D} \Phi_\rho(f_n(x)) \le c_\rho + (2-c_\rho')(\rho-\rho_n).$
\end{itemize}
Let us now define
$$G_{\rho, \alpha}= \{ u \in S_{\mu}\,|\, A(u) \leq A_{\rho}+1, B(u) \leq B_{\rho}+1 \, \mbox{and} \, |\Phi_{\rho}(u) - c_{\rho}| \leq \alpha\}.$$
We claim that, for all $\alpha >0$
$$ \inf \{  ||\Phi'_\rho|_{\M}(u)||_* \, \,|\, u \in G_{\rho, \alpha}\} =0.$$

If this is true then we obtain (i)-(iii).  This claim is proved in a standard way by using, on the special paths obtained above, a deformation argument, such as \cite[Lemma 5.15]{Willem}. We refer to \cite[Proposition 2.2]{J-PRSE1999} or \cite[Lemma 4.4]{BeRu} for details in the case of an unconstrained functional.  \smallskip

If we assume in addition that any set of the form $G:= \{ u \in \M \,|\,  A(u) \leq K_1, B(u) \leq K_2 \}$ is bounded then $(u_n) \subset E$ is bounded. Finally, under the above regularity assumptions on $\Phi_\rho: E \to \mathbb{R}$, we can pursue the proof of Theorem \ref{thm: application} starting from Step 2 and obtain $(iv)$.
\end{proof}

When the sequence $\{u_n\} \subset \M$ provided by Theorem \ref{thm: application} or Theorem \ref{thm: application-M} converges, its limit is a critical point of $\Phi_\rho|_{\M}$. In this case, the information about the approximate Morse index of $\{u_n\} \subset \M$, can be used to infer information on the Morse index of the critical point. \medskip

We recall, see for example \cite[Definition 2.5]{AcWe19}, that
\begin{definition}\label{def-index-morse} 
If $u \in \M$ is a critical point of $\Phi_\rho|_{\M}$ with Lagrange parameter $\lambda \in \R$.
\begin{itemize}
\item[1)] The Morse index $m(u) \in \N \cup \{0, \infty\}$ of $u$ with respect to $\Phi_\rho|_{\M}$ is defined as
\[
 m(u) = \sup \left\{\dim\,L\left| \begin{array}{l} \ L \text{ is a subspace of $T_u \M$ such that:~$\forall  \varphi \in L \backslash \{0\}, \,$ }
D^2\Phi_{\rho}(u) [\varphi, \varphi] < 0  \end{array}\right.\right\}.
\]

\item[2)] The free Morse index $m_f(u) \in \N \cup \{0, \infty\}$ of $u$  is defined as
\[
m_f(u) = \sup \left\{\dim\,L\left| \begin{array}{l} \ L \text{ is a subspace of $E$ such that:~$\forall  \varphi \in L \backslash \{0\}, \,$ } 
D^2\Phi_{\rho}(u) [\varphi, \varphi] < 0 \end{array}\right.\right\}.
\]
\end{itemize}
\end{definition}

\begin{remark}
Recall that, for critical points $u\in \M$ of $\Phi_\rho|_{\M}$, the restriction of $D^2\Phi_\rho(u)$ to $T_u\M$ coincides with the Hessian of $\Phi_\rho|_{\M}$ at $u$. See Remark~\ref{rema:hessian}.
\end{remark}

\begin{theorem}\label{thm: limit}
In the setting of Theorem \ref{thm: application} or of Theorem \ref{thm: application-M}, we assume that the sequence $\{u_n\} \subset \M$ converges to some $u \in \M$. Then
\begin{itemize}
\item[1)] Setting $\lambda_{\rho}:= \frac{1}{\mu}\Phi'_\rho(u)\cdot u $ 
we have
\begin{equation}\label{E-L1}
 \Phi'_{\rho}(u) - \lambda_{\rho}\langle Gu, \cdot \rangle = 0, \quad \mbox{in } E'.
\end{equation}
Here $G : E \rightarrow E$ is the injective linear map defined in \eqref{definition-G}. Equation \eqref{E-L1} indicates that
$u\in \M$ is a constrained critical point of $\Phi_{\rho}$ with Lagrange parameter $\lambda_{\rho} \in \R$. Alternatively, $u \in \M$ can be viewed as a free critical point of the functional $\Phi_{\rho}- \lambda_{\rho}Gu$ defined on $E$. 
\item[2)] The Morse index $m(u)$ of $u $ is at most $d$.
\item[3)] The free Morse index $m_f(u)$ of $u$ is at most $d+1$.
\end{itemize}
\end{theorem}
\begin{proof}
Using Equation \eqref{free-gradient} in Remark \ref{Gradient}, we immediately deduce that 1) holds.  To show that $m(u) \leq d$ we assume by contradiction that there exists a $W_0 \subset T_{u}\M$ with $\dim W_0 =d+1$ such that
$$D^2 \Phi_{\rho}(u)[w,w]  < 0, \quad \mbox{for all } w \in W_0 \backslash \{0\}. $$
Since $W_0$ is of finite dimension, by compactness and homogeneity, there exists a $\beta_0 >0$ such that
\begin{equation*}
D^2 \Phi_{\rho}(u)[w,w] < - \beta_0 ||w||^2,  \quad \mbox{for all } w \in W_0 \backslash \{0\}.
\end{equation*}
Now, from Corollary \ref{coro:step1} we deduce, for $\delta_1 >0$ given by Equation \eqref{def:delta1}, that for any 
$v\in B(\M; u,\delta_1)$,
\begin{equation}\label{L-conditionD2step1} D^2 \Phi_{\rho}(v)[w,w] < - \frac{3\beta_0}{4} ||w||^2, \quad \mbox{for all } w \in W_0 \backslash \{0\}.\end{equation}
Since $\{u_n\} \subset \M $ converges to $u$ we have that $u_n \in  B(\M; u,\delta_1)$ for $n \in \N$ large enough. Then since $\dim W_0 >d$, \eqref{L-conditionD2step1} provides a contradiction with Theorem \ref{thm: application} (iv) where we recall that $\zeta_n \to 0^+$. This proves 2). Finally, recording that $\M$ is of codimension 1 in $E$ we immediately obtain that $m_f(u) \leq d+1$.  
\end{proof}

\begin{remark}\label{extension}
If in Theorem \ref{thm: application} the conclusion only holds for almost every $\rho \in I$, this is due to the fact that it is not known if, for a given $\rho \in I$, the functional $\Phi_{\rho}$  admits a sequence of  $\{A_n\} \subset \mathcal{F}$ as in Step 1 in the proof of Theorem \ref{thm: application}. For a functional for which this is known a priori  Theorem \ref{M-theorem} directly implies the existence of sequences
$\{u_n\} \subset \M$ and $\zeta_n \to 0^+$ such that, as $n \to + \infty$, the properties (i)-(iv) in Theorem \ref{thm: application} hold.
\end{remark}

\section{The case of $\Zd$-homotopic families} \label{Theorem-Symmetric}

In this final section, we prove Theorem~\ref{Objective1}. For this we shall extend Theorem~\ref{M-theorem} to a symmetric setting. We consider the action of $\Zd$ on $\R^n$ determined by an \emph{isometric} involution of $\R^n$ with its usual distance, that we denote by $\inv$. For any subset $D\subset \R^n$ we denote by:
\[ SD=\{ Sx,\,\, x\in D\}.\]
A subset $D$ is invariant or stable if $SD=D$, in this case, a continuous map $f: D \rightarrow \M$ is said to be equivariant if:
\[f\circ S=-f.\]
We assume throughout this section that $D_0 \subset D \subset \R^n$ are compact sets such that
\begin{equation}\label{conditions}
\{ x\in D, \,\, Sx=x\}= \emptyset, \quad SD = D, \quad SD_0 = D_0.
\end{equation}

We shall restrict our attention to the following class of homotopic families:
\begin{definition}\label{C-def: def1.5 equi}
A family $\mathcal{F}$ of subsets of $\M$ will be said to be a {\it $\Zd$-homotopic family of dimension at most $n$ with boundary $B$} if there exist compact sets $D_0 \subset D \subset \R^n$ satisfying \eqref{conditions} for some isometric involution $S$ and  a continuous equivariant map $\eta_0: D_0 \rightarrow B$ such that:
\[\mathcal{F}= \{ A \subset \M \, |\, A=f(D) \mbox{ for some } f\in C(D;\M) \mbox{ with } f\circ \inv =-f \textrm{ and }  f = \eta_0 \textrm{ on $D_0$}\}.\]
\end{definition}

Our present goal in this section is to prove the following symmetric version of Theorem~\ref{M-theorem}.
\begin{theorem}\label{M-theoremSym}
Let $\phi$ be a $C^2$-functional on $E$, satisfying \ref{hol2} for some $\alpha \in (0,1]$, and assume that $\phi |_{\M}$ is even. \\
Let $\mathcal{F}$ be a $\mathbb{Z}_2$-homotopic family of $\M$ of dimension at most $n$ with boundary $B$ such that $$c:= \inf_{A \in \mathcal{F}} \max_{u \in A}\phi(u) > \max_B \phi$$ is finite.
Let $R>1$, $0< \alpha_1 \leq \frac{\alpha}{2(\alpha +2)}<1$ and $\va >0$.
Then for any $A \in \mathcal{F}$ with $\max_{u \in A} \phi(u) \leq c + \varepsilon$ satisfying 
\begin{equation}\label{top-boundedSym}
K:= \{ u \in A \,|\, \phi(u) \geq c - \va \} \subset B(0,R-1), 
\end{equation}
assuming that $\varepsilon >0$ is sufficiently small there exists $u_{\varepsilon} \in \M$ such that
\begin{enumerate}
\item $c- \va \leq \phi(u_{\va}) \leq c + \va$ ;
\item $||\phi'|_{\M}(u_{\va})|| \leq 3 \va^{\alpha_1}$ ;
\item $ u_{\va} \in A$ ;
\item If $D^2\phi(u_{\va}) [w, w]  < -  \va^{\alpha_1} \|w \|^2$ for all $w \neq 0$ in a subspace $W$ of $T_{u_{\varepsilon}}\M$, then $\dim \, W \leq n$.
\end{enumerate}
By symmetry this also holds for $-u_\varepsilon$.
\end{theorem}

The main difficulty in the proof of Theorem~\ref{M-theoremSym} is to adapt the deformation process so that $\Zd$-homotopic family $\mathcal{F}$  remains stable under our deformations. Our strategy, given $A=f(D) \in \mathcal{F}$, is to deform simultaneously near any $f(x)=u$ and $f(Sx)=-u$, transporting the data defining the deformation of $f$ near $f(x)$ to $f(Sx)$ via the antipodal map. In order to avoid any overlap when deforming, we crucially and repeatedly use the following lemma:
\begin{lemma}\label{ElementaryFact}
For any $x\in \M$, $0<\delta < \sqrt{2\mu}$ we have:
\[\emptyset=B(\M; x,\delta)\cap B(\M;-x,\delta)= B(\M; x,\delta)\cap -B(\M; x,\delta). \]
\end{lemma}

\subsection{Adapting the deformations and proof of Theorem~\ref{M-theoremSym}}
First we prove a technical  lemma that extends Lemma~\ref{lemm:3.4-FG} to the symmetric case; we recall that we assume that $\phi |_{\M}$ is even.
\begin{lemma}\label{lemm:3.4-FGSym}
Assume that the hypotheses of Lemma~\ref{lemm:3.3} are satisfied and let $f$ be a continuous equivariant map from $D$ into $\M$. Suppose that $K_1 \subset D$ is a compact subset of $D$ such that $f(K_1)\subset B(\M;u_0,\frac{\delta_4}{2})$, where $\delta_4=\min(\delta_3, 2\sqrt{2}\sqrt{\mu})$ and $\delta_3$ satisfies Eq.~\eqref{def:delta3}.  Then for sufficiently small $\nu>0$ and $t_0\in (0,t_{max})$  there is a continuous map $\eta: [0,t_{max}] \times D \rightarrow \M$ that satisfies:
\begin{enumerate}
\item $\eta(t,x)= f(x)$, if $(t,x) \in (\{0\}\times D) \cup ([0,t_{max}]\times D\setminus (N_\nu(K_1) \cup N_\nu(\inv K_1))$;
\item $\phi(\eta(t,x)) \leq \phi(f(x))$ if $(t,x) \in [0,t_{max})\times D$;
\item $\phi(\eta(t,x)) < \phi(f(x)) -\frac{\beta}{24}t^2$ if $(t,x)\in[t_0,t_{max})\times (K_1 \cup \inv K_1)$;
\item $||\eta(t,x)-f(x)|| \leq 3t $ for all $(t,x)\in [0,t_{max}]\times D$;
\item $ \eta (t,x) = - \eta(t,\inv x)$ for all $(t,x) \in [0,t_{max}]\times D$.
\end{enumerate}
\end{lemma}
\begin{proof}
Firstly, since $f$ is equivariant it follows from Lemma~\ref{ElementaryFact} that $K_1\cap \inv K_1 = \emptyset$.
Let $N_\nu(K_1)$ denote the $\nu$-neighbourhood of $K_1$ in $\mathbb{R}^n$.

Assume that $\nu >0$ is small enough so that, $f(\overline{N_\nu}(K_1)\cap D)\subset B(\M; u_0,\frac{\delta_4}{2})$ and $\overline{N_\nu}(K_1) \cap  \inv \overline{N_\nu}(K_1)\cap D = \emptyset$. Define:
\[T=\{x \in \overline{N_\nu}(K_1)\cap D \,|\,  \hat{w}(f(x))=0\}.\]
By the last point in Lemma~\ref{lemm:3.3}, for every $y\in T$ one can find $\nu^y >0$ such that for any $z\in B(y,\nu^y)\cap D$, and any $w\in T_{f(z)}W$, $||w||=1$ we have the inequality:
\[\phi(\exp_{f(z)}(tw)) < \phi(f(z)) - \frac{\beta}{24}t^2,\quad  t\in [t_0, t_{max}).\]
Put $O=\cup_{y\in T} B(y,\frac{\nu^y}{2})$ and let $h:D \rightarrow [0,1]$ be a continuous function such that: 
\[ h(x)=\begin{cases} 1 & x\in K_1 ,\\ 0 & x\in D\setminus N_\nu (K_1), \end{cases}\]
and set: $g(x)=h(x)+h(\inv x)$. Then $g$ is a continuous function on $D$ such that $g\circ \inv =g$, and, since $N_\nu(K_1)\cap \inv N_\nu(K_1) =\emptyset$, $g$ satisfies:
\[ g(x)=\begin{cases} 1 & x\in K_1\cup \inv K_1 ,\\ 0 & x\in D\setminus (N_\nu(K_1)\cup \inv N_\nu(K_1)). \end{cases}\]

\noindent
Note that since $\inv N_\nu(K_1)=N_\nu(\inv K_1)$ then $N_\nu(K_1)\cup \inv N_\nu (K_1)$ is symmetric, this follows from the fact that $\inv$ is an isometry. 

Next, proceeding as in the proof of Lemma~\ref{lemm:3.4-FG}, we choose an orthonormal basis $(e_1,\dots,e_{n+1})$ of $W \subset T_{u_0}\M$ and let $f_1: (\overline{N_\nu}(K_1) \setminus O) \cap D \rightarrow S^{n} \subset \mathbb{R}^{n+1}$ be defined by:
\[ f_1(x)=(\langle T_{f(x)} e_1, w(f(x))\rangle, \dots, \langle T_{f(x)} e_{n+1}, w(f(x))\rangle).   \]

As in the proof of Lemma~\ref{lemm:3.3}, the continuity of the map: $u\mapsto T_u|_W$ guarantees that $f_1$ is itself continuous.  
Since $(\overline{N_\nu}(K_1) \setminus O) \cap D$ is closed in $\mathbb{R}^n$ one can extend $f_1$ to a continuous map $f_2: \overline{N_\nu}(K_1)\cap D \rightarrow S^n$ by Lemma~\ref{Corollary 2.4-FG}.  Now we define:
$f_3 : (\overline{N_\nu}(K_1) \cap D) \cup  (\inv\overline{N_\nu}(K_1) \cap D) \rightarrow S^n$ by 

\[ f_3(x)=\begin{cases} f_2(x)  & x\in \overline{N_\nu(K_1)} \cap D ,\\ -f_2(\inv x) & x\in \inv\overline{N_\nu(K_1)} \cap D. \end{cases}\]

Then using Lemma~\ref{Corollary 2.4-FG} again we extend $f_3$ to a continuous map $f_4 : \R^n \rightarrow S^n$. Finally, define  a continuous map $f_5$ on $D$:
\[f_5(x) = \sum_{i=1}^{n+1} E_i^*(f_4(x))T_{f(x)}e_i, \in T_{f(x)}W, \]
where we denote $(E_i^*)$ the dual basis of the canonical basis $(E_1,\dots, E_{n+1})$ of $\mathbb{R}^{n+1}$. Now setting \[\eta(t,x)=\exp_{f(x)}(tg(x)f_5(x)), x \in D, t \in [0,t_{max}].\]
Since, $\exp_{-x}(-tv)=-\exp_x(tv)$, it follows that $\eta$ satisfies the required conditions. The point 4 follows from Equation~\eqref{eq:diffv0vf} and the mean value theorem.
\end{proof}

Next we extend Lemma~\ref{Lemma3.5-FG}:

\begin{lemma}\label{Lemma3.5-FGSym}
Let $\phi: E \rightarrow \mathbb{R}$ be a $C^2$-functional such that $\phi', \phi''$ are $\alpha$-Hölder continuous ($\alpha \in (0,1]$) on bounded sets and $\phi|_{\M}$ is even; let $f: D \rightarrow \M$ be a continuous equivariant map.
Suppose $K_2$ is a symmetric, i.e. $K_2=\inv K_2$, compact subset of $D$ with the following properties:  
\begin{itemize}
\item There exists $R >1$ such that $f(K_2) \subset B(0, R-1)\cap \M$.
\item There exists a constant $\beta >0$ such, that for all $y \in K_2$, there is a subspace $W_y$, of $T_{f(y)}\M$ with $\dim \, W_y \geq n+1$ so that 
\begin{equation}\label{3.2-FG}
D^2 \phi(f(y))[w,w]< - \beta ||w||^2, \quad \mbox{for all } \, w \in W_y \backslash \{0 \}.
\end{equation}
\end{itemize}

Then for any $0 < \delta \leq \delta_4$ where $\delta_4 >0$ was defined in Lemma~\ref{lemm:3.4-FGSym} and $\nu >0$ there is a continuous equivariant map  $\hat{f} : D \rightarrow \M$ such that if $N:= N(n)$ is the number given in 
Lemma \ref{Lemma 2.5-FG}, we have : \smallskip

\begin{itemize}
\item[(i)] $\hat{f}(x)= f(x)$ for $u \in D \backslash N_{\nu}(K_2)$; \smallskip
\item[(ii)] $\phi(\hat{f}(x)) \leq \phi (f(x))$  for all $x \in D$; \smallskip
\item[(iii)] If $x \in K_2$, then $\phi(\hat{f}(x)) <  \displaystyle \phi (f(x)) - \frac{\beta \, \delta^2}{3456N^2};$ \smallskip
\item[(iv)] $ ||\hat{f}(x) - f(x) || \leq \frac{\delta}{2}$ for all $x \in D$.
\end{itemize}
\end{lemma}

\begin{proof}
Let $0 < \delta \leq \delta_4$ be fixed. 
As in the proof of Lemma~\ref{Lemma3.5-FG}, we can find a $\varepsilon_N >0$ and points $x_1, \cdots, x_m$ in $\R^n$ such that $\cup_{i=1}^{m}B(x_i, \varepsilon_N)$ covers $K_2$ and $B(x_i,\varepsilon_N)\cap K_2\neq\emptyset$ for each $i \in\{1, \cdots, m\}$. Choosing $y_i \in B(x_i, \varepsilon_N) \cap K_2$ and setting $B_{y_i} = B(x_i, \varepsilon_N)$ we observe, as before, that, taking $ \nu >0$ small enough, we can assume that $\cup_{i=1}^m N_{\tau}(B_{y_i}) \subset N_{\nu}(K_2)$,  and for $1 \leq i \leq m,$
\begin{equation}\label{equi}
f(x) \in B \Big(\M; f(y_i), \frac{\delta}{4N} \Big), \quad x \in \overline{N}_{\tau}(B_{y_i}) \cap D.
\end{equation}
In addition, any intersection of $N$ distinct sets $\overline{N_{\tau}}(B_{y_i})$ is empty. Note that so far we have not used the assumption that $K_2$ is symmetric.

We shall now define by induction, continuous equivariant functions $f_0, f_1, \cdots, f_m : D \rightarrow \M$ such that for all $1 \leq i \leq m$ we have that
\begin{equation}\label{3.3}
\phi(f_i(x)) < \phi(f_{i-1}(x)) - \frac{\beta \, \delta^2}{3456 N^2} \quad \mbox{if } x \in (\overline{B}_{y_{k+1}}\cup \inv\overline{B}_{y_{k+1}}) \cap D,
\end{equation}
\begin{equation}\label{3.4}
\phi (f_i(x)) \leq \phi (f_{i-1}(x)) \quad \mbox{if } x \in  D,
\end{equation}
and
\begin{equation}\label{3.5}
||f_i(x) - f_{i-1}(x)|| \leq
\begin{cases}
0  \quad \mbox{if } x \in D \backslash (N_{\tau}(\overline{B}_{y_i})\cup \inv N_{\tau}(\overline{B}_{y_i})),\\
\frac{\delta}{4 N}\quad \mbox{if } x \in (N_{\tau}(\overline{B}_{y_i})\cup \inv N_{\tau}(\overline{B}_{y_i})) \cap D.
\end{cases}
\end{equation}
Let $f_0= f$ and suppose that $f_0, f_1, \cdots, f_k$ are well-defined and satisfy inequalities \eqref{3.3}, \eqref{3.4} and \eqref{3.5} for $k <m.$ Clearly
$$||f_i(x) - f(x) || \leq \frac{i \delta}{4 N} \quad \mbox{if } x \in \bigcap_{j=1}^i (N_{\tau}(\overline{B}_{y_i})\cup \inv N_{\tau}(\overline{B}_{y_i})) \cap D.$$   \medskip

Since any intersection of $N$ distinct $N_{\tau}(\overline{B}_{y_i})$ is empty, by symmetry the same is true for $\inv N_{\tau}(\overline{B}_{y_i})$ and hence the above intersection is empty whenever $i> 2(N-1).$
Thus:
$$ ||f_k(x) - f(x) || \leq \frac{\delta(2(N-1))}{4  N} \quad \mbox{if } x \in  D.$$

%
As $f : \overline{B}_{y_{k+1}} \cap D \rightarrow B(\M; f(y_{k+1}), \frac{\delta}{4  N})$, we see that $f_k$ maps $\overline{B}_{y_{k+1}} \cap D$ into 
$B(\M; f(y_{k+1}), \frac{\delta}{2 }(1 - \frac{1}{N}+ \frac{1}{2 N}) \subset B(\M; f(y_{k+1}), \frac{\delta}{2 }).$ 
By assumption \eqref{3.2-FG}, there is some subspace $W_{y_{k+1}}$ of $E$ with $\dim W_{y_{k+1}} \geq n+1$ such that for any $w \in W_{y_{k+1}},$ with $||w|| =1,$ we have that 
$D^2 \phi(f(y_{k+1}))[w,w]< - \beta$.
Hence, we may apply Lemma \ref{lemm:3.4-FGSym} to $f_k$ and any $t_0 \in (0, t_{max})$ to obtain a continuous equivariant deformation $\eta(t,x)$ satisfying the conclusion of that lemma. Define now $f_{k+1}(x) = \eta (\frac{\delta}{12 N},x)$ to get a continuous function  $f_{k+1}: D \rightarrow \M$ satisfying
$$ \phi(f_{k+1}(x)) < \phi(f_{k}(x)) - \frac{\beta \, \delta^2}{3456 N^2} \quad \mbox{for } x \in (\overline{B}_{y_{k+1}}\cup \inv \overline{B}_{y_{k+1}}) \cap D,$$
\begin{equation*}
\phi (f_{k+1}(x)) \leq \phi (f(x)) \quad \mbox{for } x \in  D,
\end{equation*}
and
\begin{equation*}
||f_{k+1}(x) - f_{k}(x)|| \leq
\begin{cases}
0  \quad \mbox{if } x \in D \backslash (N_{\tau}(\overline{B}_{y_{k+1}})\cup  \inv N_\tau(\overline{B}_{y_{k+1}})),\\
\frac{\delta}{4 N}\quad \mbox{if } x \in (N_{\tau}(\overline{B}_{y_{k+1}})\cup \inv N_{\tau}(\overline{B}_{y_{k+1}})) \cap D.
\end{cases}
\end{equation*}
By induction we see that $f_0, \cdots, f_m$ are well-defined. Clearly $\hat{f} = f_m$ verifies the claims of the lemma.
\end{proof}

Finally, we shall need ~\cite[Lemma 3.1]{Willem}:
\begin{lemma}\label{Lemma3.6-FGSym}
Let $\phi$ be a $C^1$ functional on $E$ and let $f: D \rightarrow \M$ be a continuous equivariant map.  Let $\tilde{c}$, $\tilde{\va}$, $\tilde{\mu} >0$
 be three constants.
 Suppose $K_3$ is a symmetric compact subset of $D$ satisfying
$$ \tilde{c} - \tilde{\va} \leq \phi(f(x)) \leq \tilde{c} + \tilde{\va}, \quad \mbox{for } x \in K_3.$$
Assume that, for all $x \in K_3$,
$$||\phi'|_{\M}(u) || \geq \frac{8 \tilde{\va}}{\tilde{\mu}}, \quad \mbox{for } u \in B(\M; f(x), 2 \tilde{\mu}), $$
then there is a equivariant continuous map  $\hat{f} : D \rightarrow \M$ such that 
\smallskip

\begin{itemize}
\item[(i)] $\hat{f}(x)= f(x)$ if $ \,\phi(f(x)) \leq \tilde{c} - 2 \tilde{\va}$; \smallskip
\item[(ii)] $\phi(\hat{f}(x)) \leq \phi (f(x))$  for all $x \in D$; \smallskip
\item[(iii)] If $x \in K_3$, then $\phi(\hat{f}(x))\leq  \tilde{c} - \tilde{\va}.$
\end{itemize}
\end{lemma}

\begin{proof}[Proof of Theorem~\ref{M-theoremSym}]
Using the notations introduced in the proof of Theorem~\ref{M-theorem} we observe that the sets:
\[
\begin{gathered}
K = \{ x \in D \,|\, \phi(f(x)) \geq c- \varepsilon \},\\
T_1= \{x \in K \,|\, ||\phi'|_{\M}(u)|| > \varepsilon^{\alpha_1}, \mbox{ for all } u \in B(\M;f(x),\hat{\delta})\}, \\
T_2=K\setminus T_1,
\end{gathered}
\]
are symmetric since $\phi$ is even. Note also that for any symmetric subset $A \subset D$ then $\overline{A}$ is symmetric.

The arguments of the proof are then identical to those of the proof of Theorem~\ref{M-theorem}, using instead Lemmata~\ref{lemm:3.4-FGSym}, \ref{Lemma3.5-FGSym} and~\ref{Lemma3.6-FGSym} and modifying constants where necessary.
\end{proof}

\begin{remark}
Inspecting the proofs of Lemmata~\ref{lemm:3.4-FGSym}, \ref{Lemma3.5-FGSym} and~\ref{Lemma3.6-FGSym}, one might observe that the assumption that $S$ is an isometry can be omitted if we replace the usual distance $d$ of $\R^n$ by a topologically equivalent distance for which $S$ is an isometry and such that the conclusions of Lemma~\ref{Lemma 2.5-FG} hold. This would be the case of the \enquote{average} distance:
\[ d_*(x,y)= \frac{1}{2}\left(d(x,y)+d(Sx,Sy)\right). \]
\end{remark}

\begin{proof}[Proof of Theorem~\ref{Objective1}]
The proof follows exactly that of Theorem \ref{thm: application}, with Theorem \ref{M-theoremSym} replacing Theorem~\ref{M-theorem}. We apply Theorem \ref{M-theoremSym} with the particular choices: \[\inv\cdot(s,x)=(s,-x), \quad (s,x) \in \R \times \R^{n-1}; \quad D=[0,1]\times S^{N-2} \subset \R^N, D_0= \{0,1\}\times S^{N-2}.\] 
Note that $\gamma_0, \gamma_1 : S^{N-2} \rightarrow S_{\mu}$ corresponds to $\eta_0$.
\end{proof}

\bibliographystyle{plain}
\bibliography{biblio}
\vspace{0.25cm}

\end{document}